\documentclass{amsart}

\usepackage{amssymb}
\usepackage{amsfonts}
\usepackage{amsmath}
\usepackage{graphicx}%
\usepackage{lastpage}
\usepackage[utf8]{inputenc}
\usepackage[legalpaper,bookmarks=true,colorlinks=true,linkcolor=blue,citecolor=blue]{hyperref}
\usepackage{fancyhdr}
\usepackage{color}
\usepackage[mathlines]{lineno}
\usepackage{multirow}
\usepackage{caption}
\usepackage{subcaption}
\usepackage{lscape}
\usepackage{epsfig}
\usepackage{natbib}
\usepackage{float}

\newtheorem{theorem}{Theorem}
\theoremstyle{plain}

\newtheorem{corollary}{Corollary}

\newtheorem{definition}{Definition}

\numberwithin{equation}{section}

\begin{document}

\title{Divergence Measures Estimation and Its Asymptotic Normality Theory in the discrete case}
\author{$^{(1)}$ BA Amadou Diadi\'e}
\email{ba.amadou-diadie@ugb.edu.sn}
\author{$^{(1,2,4)}$LO Gane Samb}
\email{gane-samb.lo@ugb.edu.sn}
\begin{abstract}
In this paper we provide the asymptotic theory
of the general of $\phi$-divergences measures, which includes the most common divergence measures: Renyi and Tsallis families and the Kullback-Leibler measure. We are interested in divergence measures in the discrete case. One sided and two-sided statistical tests are derived as well as symmetrized estimators. Almost sure rates of convergence and asymptotic normality theorem are obtained in the general case, and next particularized for the Renyi and Tsallis families and for the Kullback-Leibler measure as well. Our theorical results are validated by simulations.  \end{abstract}

\maketitle
\section{Introduction}
\subsection{Motivations}
$ $ \\

\noindent In this paper, we study the convergence of $\phi-$divergence measure estimator for empirical discrete probability distributions supported on a finite set.\\

\noindent Let throughout the following $\mathcal{X}=\{c_1,c_2,\cdots,c_r\}\,(r\geq 2)$ be a finite countable space. The probability distributions on $\mathcal{X}$ are finite dimensional vectors $\textbf{p}$ in $$\mathcal{P}(\mathcal{X})=\left\lbrace \textbf{p}=(p_c)_{c\in \mathcal{X}}:p_c\geq 0,\,\forall c\in \mathcal{X}\ \ \text{and}\ \ \sum_{c\in \mathcal{X}}p_c=1\right\rbrace.$$
A divergence measure on $\mathcal{P}(\mathcal{X})$\ is a function 
\begin{equation}
\begin{tabular}{cccl}
$\mathcal{D}:$ & $(\mathcal{P}(\mathcal{X}))^{2}$ & $\longrightarrow $ & $\overline{\mathbb{%
R}}$ \\ 
& $(\textbf{p},\textbf{q})$ & $\longmapsto $ & $\mathcal{D}(\textbf{p},\textbf{q})$
\end{tabular}
\label{divApp}
\end{equation}

\noindent such that $\mathcal{D}(\textbf{p},\textbf{p})=0$ for any $\textbf{p}$ such that $(\textbf{p},\textbf{p})$ in the domain of application of 
$\mathcal{D}$.\\

\noindent The function $\mathcal{D}$ is not necessarily a mapping. And if it is, it is not always symmetrical and it does neither have to be a metric. In lack of symmetry, the following more general notation is more appropriate :

\begin{equation}
\begin{tabular}{cccl}
$\mathcal{D}:$ & $\mathcal{P}_1(\mathcal{X}) \times \mathcal{P}_2(\mathcal{X})$ & $\longrightarrow $ & $\overline{\mathbb{%
R}}$ \\ 
& $(\textbf{p},\textbf{q})$ & $\longmapsto $ & $\mathcal{D}(\textbf{p},%
\textbf{q})$,
\end{tabular}
\label{divAppG}
\end{equation}

\noindent where $\mathcal{P}_1(\mathcal{X})$ and $\mathcal{P}_2(\mathcal{X})$ are two families of probability distributions  on 
$\mathcal{X} $, not necessarily the same. To better explain our concern, let us introduce some of
the most celebrated divergence measures.\\

\bigskip \noindent Let $(\textbf{p},\textbf{q})\in \mathcal{P}(\mathcal{X})\times \mathcal{P}(\mathcal{X})$ with $\mathcal{X}=\{c_1,c_2,\cdots,c_r\}$,  and let $X$ and $Y$ two randoms variables such that 
$$\mathbb{P}(X=c_j)=p_j,\ \ \text{and}\ \ \mathbb{P}(Y=c_j)=q_j,\ \  j\in \{1,\cdots,r\}.$$

\noindent and set $\textbf{p}=(p_1,\cdots,p_r)^t$ and $\textbf{q} =(q_1,\cdots,q_r)^t$.
\\

\bigskip \noindent The four most popular divergence are :
\\

\noindent (1) The $L_{2}^{2}$-divergence measure :
\begin{equation}
\mathcal{D}_{L_{2}}(\textbf{p},\textbf{q})=\sum_{j=1}^r(p_j-q_j)^2 . \label{L22}
\end{equation}

\bigskip \noindent (2) The family of Renyi's divergence measures indexed by $\alpha \neq 1$, $\alpha>0$, known under the name of Renyi-$\alpha$ : 
\begin{equation}
\mathcal{D}_{R,\alpha }(\textbf{p},\textbf{q})=\frac{1}{\alpha -1}\log
\left( \sum_{j=1}^rp_j^\alpha q_j^{1-\alpha }\right). \label{renyi}
\end{equation}

\bigskip \noindent (3) The family of Tsallis divergence measures indexed by $\alpha \neq 1$, $\alpha>0$, also known under the name of Tsallis-$\alpha$ : 
\begin{equation}
\mathcal{D}_{T,\alpha }(\textbf{p},\textbf{q})=\frac{1}{\alpha -1}\left(\sum_{j=1}^rp_j^\alpha q_j^{1-\alpha }-1\right) ;  \label{tsal}
\end{equation}

\bigskip \noindent (4) The Kulback-Leibler divergence measure
\begin{equation}
\mathcal{D}_{KL}(\textbf{p},\textbf{q})=\sum_{j=1}^rp_j\log(p_j/q_j).
\label{kull1}
\end{equation}

\bigskip \noindent The latter, the Kullback-Leibler divergence measure, may be interpreted as a limit case
of both the Renyi's family and the Tsallis' one by letting $\alpha
\rightarrow 1$. As well, for $\alpha $ near 1, the Tsallis family may be
seen as derived from $\mathcal{D}_{R,\alpha }(\textbf{p},\textbf{q})$ based on the first order expansion of the logarithm
function in the neighborhood of the unity. Here for ease of notation we refer the notation log as the natural logarithm \newline

\noindent From this small sample of divergence measures, we may give the following remarks :\\

\noindent For both the Renyi and the Tsallis families, we may have computation problems. So without loss of generality, suppose
\begin{equation}
\ \ \ p_j>0\ \ \text{and}\ \ q_j>0, \ \ \forall j\in D=\{1,2,\cdots,r\}\ \ \ \ \ (\textbf{BD})\label{BD}
\end{equation}

\bigskip \noindent If Assumption (\ref{BD}) holds, we do not have to worry about summation problems, especially for Tsallis, Renyi and Kulback-Leibler measures, in the computations arising in estimation theories. 
 This explains why Assumption (\ref{BD}) is systematically used in a great number of works in that topic, for example, in \cite{singh}, \cite{kris}, \cite{hall}, and recently in \cite{ba1} to cite a few.\\

\noindent \noindent It is clear from the very form of these divergence measures that we do not have symmetry, unless for the special case where $\alpha=1/2$. So we define the following symetric version of divergence measures
$$\mathcal{D}^{(s)}(\textbf{p},\textbf{q})=\frac{\mathcal{D}(\textbf{p},\textbf{q})+\mathcal{D}(\textbf{q},\textbf{p})}{2}$$
provided that $\mathcal{D}(\textbf{p},\textbf{q})$ and $\mathcal{D}(\textbf{q},\textbf{p})$ are finite.
\\
 
\bigskip \noindent Both families are build on the following summation 
\begin{equation*}
\mathcal{S}_{\alpha}(\textbf{p},\textbf{q})=\sum_{j\in D}p_j^{\alpha }q_j^{1-\alpha },\ \ \text{with}\ \ \alpha \neq 1,\ \ \alpha>0.
\end{equation*}

\noindent Although we are focusing on the aforementioned divergence measures in this paper, it is worth mentioning that there exist quite a few number of them. Let us cite for example the ones named after  : Ali-Silvey or $f$-divergence \cite{topsoe}, Cauchy-Schwarz, Jeffrey divergence (see \cite{evren}), Chernoff (See \cite{evren}) , Jensen-Shannon (See \cite{evren}). According to \cite{cichocki}, there is more than a dozen of different divergence measures in the literature.\\

\noindent Before coming back to our divergence measures estimation of interest, we want to highlight some important applications of them. Indeed, divergence has proven to be useful in applications. Let us cite a few of them :\\

\noindent (a) They heavily intervene in Information Theory and recently in Machine Learning.\\

\noindent (b) They have been used as similarity measures in image registration or multimedia classification (see \cite{moreno}).\\

\noindent (c) They are also used as loss functions in evaluating and optimizing the performance of density estimation
methods (see \cite{hall}).\\

\noindent (d) Divergence estimates can also be used to determine sample sizes required to achieve given performance levels in hypothesis testing.\\

\noindent (e) There has been a growing interest in applying divergence to various fields of science and engineering for the purpose of estimation,
classification, etc. (See \cite{bhattacharya}, \cite{liu1}).\\

\noindent (f) Divergence also plays a central role in the frame of large deviations results including the asymptotic rate of decrease of error probability in binary hypothesis testing problems.\newline

\noindent (g) The estimation of divergence between the samples drawn from unknown distributions gauges the distance between those distributions.
Divergence estimates can then be used in clustering and in particular for deciding whether the samples come from the same distribution by comparing
the estimate to a threshold.\\

\noindent (h) Divergence gauges how differently two random variables are distributed and it provides a useful measure of discrepancy between
distributions. In the frame of information theory, the key role of divergence is well known.\\

\noindent The reader may find more applications and descriptions in the following papers :  \cite{kullback},\cite{fukunaga}, 
\cite{cardoso}, \cite{ojala}, \cite{hastie},  \cite{moreno},\cite{macKay}.\\

\bigskip \noindent In the next subsection, we describe the frame in which we place the estimation problems we deal in this paper.\\

\subsection{Statistical Estimations} \label{subsec_intro_estim}$ $\\

\noindent The divergence measures may be applied to two statistical problems among others.\\

\noindent \textbf{(A)} First, it may be used as a fitting problem as described here. Let $X_{1},X_{2},\cdots$ a sample of replications of $X$ with an unknown probability distribution $\textbf{p}$ and we want to test the hypothesis that $\textbf{p}$ is equal to a known and fixed probability $\textbf{p}_0.$ Theoretically, we can answer this question by estimating a divergence measure $\mathcal{D}(\textbf{p},\textbf{p}_{0})$ by a plug-in estimator $\mathcal{D}(\widehat{\textbf{p}}_n,\textbf{p})$ where, for each $n\geq 1$, $\textbf{p}$ is replaced by an estimator $\widehat{\textbf{p}}_n$ of the probability law, which is based on sample $X_1$, $X_2$, ..., $X_n$, to be precised.\\

\noindent From there establishing an asymptotic theory of $\Delta _{n}=\mathcal{D}(%
\widehat{\textbf{p}}_n,\textbf{p}_0)-\mathcal{D}(\textbf{p},\textbf{p}_{0})$ is thought to be necessary to conclude.\\

\noindent \textbf{(B)} Next, it may be used as tool of comparing for two distributions. We may have two samples and wonder whether they come from the same probability distribution. Here, we also
may two different cases.\\

\noindent \textbf{(B1)} In the first, we have two independent samples $%
X_{1},X_{2},....$ and $Y_{1},Y_{2},....$ respectively from a random variable 
$X$ and $Y$ according the probability distributions $\textbf{p}$ and $\textbf{q}$. Here the estimated divergence $\mathcal{D}(\widehat{\textbf{p}}_n,%
\widehat{\textbf{q}}_m)$, where $n$ and $m$ are the sizes of the available samples, is the natural estimator of $\mathcal{D}(\textbf{p},\textbf{q})$ on which depends the statistical test of the hypothesis : $\textbf{p}=\textbf{q}$.\\

\noindent \textbf{(B2)} But the data may also be paired $(X,Y)$, $(X_{1},Y_{2}),(X_{2},Y_{2}),...,$
that is $X_{i}$ and $Y_{i}$ are measurements of the same case $i=1,2,...$ In
such a situation, testing the equality of the margins $\textbf{p}_X=\textbf{p}_Y$ should be based on an estimator $\widehat{\textbf{p}}_{X,Y}^{(n)}$ of the joint probability law of the couple $(X,Y)$ based of the
paired observations $(X_{i},Y_{i})$, $i=1,2,\ldots,n$.\\

\bigskip \noindent We did not encounter the approach (B2) in the literature. In the (B1) approach, almost all the papers used the same sample size, at the exception of \cite{poczos}, for the double-size estimation problem. In our view, the study case should rely on the available data so that using the same sample size may lead to a loss of information. To apply their method, one should take the minimum of the two sizes and then loose information. We suggest to come back to a general case and then study the asymptotic theory of $\mathcal{D}(\widehat{\textbf{p}}_n,\widehat{\textbf{q}}_m)$ based on samples $X_{1},X_{2},..,X_n.$ and $Y_{1},Y_{2},...,Y_m$. In this paper, we will systematically use arbitrary samples sizes.

\subsection{Previous work}$ $ \\

\noindent In the context of the situation (B1), there are several papers dealing with the estimation of the divergence measures. As we are concerned in this paper by the weak laws of the estimators, our review on that problematic did not return significant things. Instead, the literature presented us many kinds of results on almost-sure efficiency of the estimation, with rates of convergences and laws of the iterated logarithm, $L^{p}$ ($p=1,2$) convergences, etc. To be precise, \cite{dakher} used recent techniques based on functional empirical process to provide a series of interesting rates of convergence of the estimators in the case of one-sided approach for the class de Renyi, Tsallis, Kullback-Leibler to cite a few. Unfortunately, the authors did not address the problem of integrability, taking that the divergence measures are finite. Although the results should be correct under the boundedness assumption \eqref{BD} (\textbf{BD})  we described earlier, a new formulation in that frame would be welcome.\\

\bigskip \noindent In the context of the situation (B1), we may cite first the works of \cite{kris} and \cite{singh}. They both used divergence measures based on probability density functions and concentrated of Renyi-$\alpha $, Tsallis-$\alpha $ and Kullback-Leibler.

\bigskip \noindent Specifically, \cite{kris} defined Reyni and Tsallis estimators by correcting the plug-in estimator and established that, as long as 
$\mathcal{D}_{R,\alpha }(\textbf{p},\textbf{q}) \geq c$ and $\mathcal{D}_{T,\alpha }(\textbf{p},\textbf{q}) \geq c$, for some constant $c>0$, then
\begin{eqnarray*}
&& \mathbb{E}\left \vert \mathcal{D}_{R,\alpha} (\widehat{\textbf{p}}_n,\widehat{\textbf{q}}_n)-\mathcal{D}_{R,\alpha }(\textbf{p},\textbf{q})\right\vert  \leq  c \left( n^{-1/2}+n^{-\frac{3s}{2s+d}}
\right)\\ 
\text{and}&& \\
&& \mathbb{E}\left \vert \mathcal{D}_{T,\alpha }(\widehat{\textbf{p}}_n,\widehat{\textbf{q}}_n)-\mathcal{D}_{T,\alpha }(\textbf{p},\textbf{q})\right\vert  \leq  c \left( n^{-1/2}+n^{-\frac{3s}{2s+d}}
\right),
\end{eqnarray*}

\bigskip \noindent \cite{poczos} used a $k-$nearest-neighbor approach to prove that 
if $|\alpha -1|<k$, ($\alpha\neq 1)$ then 
\begin{eqnarray*}
&& \lim_{n,m\rightarrow \infty }\mathbb{E}\left[ \mathcal{D}_{T,\alpha }(\widehat{\textbf{p}}_n,\widehat{\textbf{q}}_m)-\mathcal{D}_{T,\alpha }(\textbf{p},\textbf{q})
\right] ^{2}=0\\
\text{and} &&\\
&&
\lim_{n,m\rightarrow \infty }\mathbb{E}\left( \mathcal{D}_{R,\alpha }(\widehat{\textbf{p}}_n,\widehat{\textbf{q}}_m)\right) =\mathcal{D}_{R,\alpha }(\textbf{p},\textbf{q}).
\end{eqnarray*}

\bigskip \noindent There has been recent interest in deriving convergence rates for divergence estimators \cite{moon}-\cite{kris}. The rates are typically
derived in terms of smoothness $s$ of the densities : 

 \bigskip \noindent  The estimator of \cite{liu2} 
converges at rate $n^{-\frac{s}{s+d}}$, achieving the parametric rate when $s>d$.

\bigskip \noindent Similarly, \cite{sricharan} show that when $s>d$ a $k$-nearest-neighbor style estimator achieves rate $n^{-2/d}$ (in absolute error)
ignoring logarithmic factors. In a follow up work, the authors improve this result to $O(n^{-1/2})$ using an
ensemble of weak estimators, but they require $s > d$ orders of smoothness.

\bigskip \noindent \cite{singh} provided an estimator for R\'enyi$-\alpha $ divergences as well as general density functionals that uses
a \textit{mirror image} kernel density estimator. They obtained exponential inequalities for the deviation of the estimators from the true value.\\

 \bigskip \noindent \cite{kall} studied an $\varepsilon-$nearest neighbor estimator for the $L_2-$divergence that enjoys the same rate of
convergence as the projection-based estimator of \cite{kris}.

\subsection{Main contributions} $ $ \\

\noindent Our main contribution may be summurized as follows, for data  sampled from one or two unknown random variables, we derive almost sure convergency and central limit theorems for empirical $\phi-$ divergences. We will focus on divergence measures between discrete probability distribution. As well, our results applied to the approaches (A) and (B1) defined above. As a consequence, we estimate divergence measures by their plug-in counterparts, meaning that we replace the probability mass function (\textit{p.m.f.}) in the expression of the divergence measure by a nonparametric estimator of the \textit{p.m.f.}'s.

\bigskip \noindent We also wish to get first general laws for an arbitrary functional of the form 

\begin{equation} \label{defJ}
J(\textbf{p},\textbf{q})=\sum_{j\in D} \phi(p_j,q_j),
\end{equation}
\noindent where $\phi : (0,1)^2\rightarrow \mathbb{R}$ is a twice continously differentiable function. The results on the functional $J(\textbf{p},\textbf{q})$, which is also known under the name of $\phi$-divergence, will lead to those on the particular cases of the Tsallis and Kullback-Leibler measures.\\

\subsection{Overview of the paper}$ $ \\

\noindent The rest of the paper is organized as follows. In \textsc{Subsection} \ref{notation_mainresults}, we define estimators of the \textit{p.m.f.} $p_j$ and $q_j$  based on i.i.d. samples according respectively to $\textbf{p}$ and $\textbf{q}$. In  \textsc{Section} \ref{mainresults}, we will give our foul results for functional $J(\textbf{p},\textbf{q})$ both one sided and two-sided approaches.
In Section \ref{particular-divergencemeasure}, we will particularize the results for specific measures we already described. \textsc{Section} \ref{proofs} provides  the proofs and in \textsc{Section} \ref{simulation} we present some simulations confirming our results. Finally 
in \textsc{Section} \ref{conclus}, we conclude.

\section{Empirical $\phi-$ divergence}
\subsection{Notations and main results}\label{notation_mainresults}$ $ \\

\noindent Before we state the main results we need a few definitions. Let $X$ and $Y$ two randoms variables defined on the probability distributions $(\mathcal{X},\mathcal{A},\mathbb{P})$ with $\mathcal{X}=\{c_1,c_2,\cdots,c_r\}$ and $\mathbf{p}=(p_j)_{1\leq j\leq r}$ and $\mathbf{q}=(q_j)_{1\leq j\leq r}$ two discrete probability distributions on $\mathcal{X}$
such that, for any $j\in D=\{1,2,\cdots,r\}$
$$p_j=\mathbb{P}(X=c_j)\ \ \text{and}\ \  q_j=\mathbb{P}(Y=c_j).$$

\noindent We suppose that \eqref{BD} is satisfied that is  $\forall\,j\in D,\ \ p_j>0\ \ \text{and}\ \ q_j>0.$\\

\noindent Define the empirical probability distribution generated by i.i.d. random variables $X_1,\cdots,X_n$ from the probability distribution $\textbf{p}$ as
\begin{eqnarray}\label{pn}
\widehat{\textbf{p}}_n=(\widehat{p}_n^{c})_{c\in \mathcal{X}},\ \ \text{where}\ \  \widehat{p}_n^{c_j}&=&\frac{1}{n}\sum_{i=1}^n1_{c_j}(X_i)
\end{eqnarray}where 
 $1_{c_j}(X_i)=\begin{cases}
 1\ \ \text{if}\ \ X_i=c_j\\
 0\ \ \text{otherwise}
 \end{cases} $ for any $j\in \{1,\cdots,r\}$. \\

\bigskip  \noindent 
$\widehat{\textbf{q}}_m$ is defined in the same way by $Y_1,\cdots,Y_m\stackrel{i.i.d.}{\sim}\textbf{q}$ that is \begin{equation}\label{qn}
\widehat{\textbf{q}}_m=(\widehat{q}_m^{c})_{c\in \mathcal{X}},\ \ \text{where}\ \  \widehat{q}_m^{c_j}=\frac{1}{m}\sum_{i=1}^m1_{c_j}(Y_i),
\end{equation}
 
\subsection{ $\phi$-divergence measure}
\label{phi_divergence}
\begin{definition}
 The $\phi$-divergence between the two probability distributions $\textbf{p}$ and $\textbf{q}$ is given by
 \begin{equation}
J(\textbf{p},\textbf{q})=\sum_{j\in D}\phi(p_j,q_j) 
 \end{equation} where $\phi:[0,1]^2\rightarrow\mathbb{R} $ is a measurable function  having continuous second order partial derivatives. 
 \end{definition}The results on the functional $J(\textbf{p},\textbf{q})$ will lead to those on the particular cases of the Renyi, Tsallis, and Kullback-Leibler measures. \\
 
\bigskip \noindent  Based on \eqref{pn} and \eqref{qn}, we will use the following empirical $\phi$-divergences. 
\begin{eqnarray*}
J(\widehat{\textbf{p}}_n,\textbf{q})&=&\sum_{j\in D}\phi(\widehat{ p}_n^{c_j},q_j), \text{\ \ \ \ }%
J(\textbf{p},\widehat{\textbf{q}}_m)=\sum_{j\in D}\phi(p_j,\widehat{q}_m^{c_j}),\\
\text{\ \ and \ \ }&& J(\widehat{\textbf{p}}_n,\widehat{\textbf{q}}_m)=\sum_{j\in D}\phi(\widehat{p}_n^{c_j},\widehat{q}_m^{c_j}).
\end{eqnarray*}

Set \begin{eqnarray}\label{abcn}
a_n&=&\sup_{j\in D}|\widehat{p}_n^{c_j}-p_j|,\ \ \ \ b_m=\sup_{j\in D}|\widehat{q}_m^{c_j}-q_j|,\\
 \ 
\notag \text{and}\ \ && c_{n,m}=\max(a_n,b_m).
\end{eqnarray}

 \noindent Denote \begin{equation*}
\phi _{1}^{(1)}(s,t)=\frac{\partial \phi }{\partial s}(s,t),\text{ }\phi
_{2}^{(1)}(s,t)=\frac{\partial \phi }{\partial t}(s,t)
\end{equation*}

\noindent and

\begin{equation*}
\phi _{1}^{(2)}(s,t)=\frac{\partial ^{2}\phi }{\partial s^{2}}(s,t),\text{ }%
\phi _{2}^{(2)}(s,t)=\frac{\partial ^{2}\phi }{\partial t^{2}}(s,t),\text{ }%
\phi _{1,2}^{(2)}(s,t)=\phi _{2,1}^{(2)}(s,t)=\frac{\partial ^{2}\phi }{%
\partial s\partial t}(s,t).
\end{equation*}
\noindent  Set \begin{eqnarray}
\label{a1pq} A_{1,p}&=&\sum_{j\in D}|\phi _{1}^{(1)}(p_j,q_j)|,\ \ \ A_{2,q}=\sum_{j\in D}|\phi _{2}^{(1)}(p_j,q_j)|\\
\label{a3pq}A_{3,q}&=&\sum_{j\in D}|\phi _{1}^{(1)}(q_j,p_j)|,\ \ \text{and}\ \  A_{4,p}=\sum_{j\in D}|\phi _{2}^{(1)}(q_j,p_j)|.
\end{eqnarray}

\section{Statements of the main results
}
 
\label{mainresults}
\subsection{Main results}$ $ \\

\noindent Here are our main results. The first concerns the almost sure efficiency of the estimators.\\

\begin{theorem} \label{thJ12} Let $\textbf{p}$ and $\textbf{q}$ two probability distributions and  
$\widehat{\textbf{p}}_n$ and $\widehat{\textbf{q}}_m$ be generated by i.i.d. samples $X_1,\cdots,X_n$ and $Y_1,\cdots,Y_m$ according respectively to $\textbf{p}$ and $\textbf{q}$ and given by \eqref{pn} and \eqref{qn}, $(BD)$ \eqref{BD} be satisfied. 
Then the following asymptotic results hold
\begin{itemize}
\item[(a)] One sample

\begin{eqnarray}
&&\limsup_{n\rightarrow +\infty }\frac{|J(\widehat{\textbf{p}}_n,\textbf{q})-J(\textbf{p},\textbf{q})|}{a_{n}}\leq A_{1,p}
,\ \ \text{a.s.}  \label{thJ12c1}
\\
&& \limsup_{m\rightarrow +\infty } \frac{\left\vert
J(\textbf{p},\widehat{\textbf{q}}_m)-J(\textbf{p},\textbf{q})\right\vert}{b_{m}} \leq A_{2,q},\ \ \text{a.s.}  \label{thJ12c2}
\end{eqnarray}
\item[(b)] Two samples :
\begin{eqnarray}
&& \limsup_{(n,m) \rightarrow (+\infty,+\infty)} \frac{\left\vert
J(\widehat{ \textbf{p}}_n,\widehat{\textbf{q}}_m)-J(\textbf{p},\textbf{q})\right\vert}{%
c_{n,m}} \leq A_{1,p}+A_{2,q} \text{\ \ a.s.}  \label{thJ22c1}
\end{eqnarray}

\end{itemize} 

\bigskip
\noindent where $a_{n}$, $b_m$ and $c_{n,m}$ are as in \eqref{abcn} and $A_{1,p}$ and $A_{2,q}$ as in \eqref{a1pq}.
\end{theorem}

\bigskip \bigskip
\noindent The second concerns the asymptotic normality of the estimators.\\

 \noindent Let 
 $$V_{1,p}=\sum_{j\in D}p_j(1-p_j)(\phi_{1}^{(1)}(p_{j},q_{j}))^2- 2\sum_{(i,j)\in D^2,\,i\neq j}p_ip_j\phi_{1}^{(1)}(p_{i},q_{i})\phi_{1}^{(1)}(p_{j},q_{j})$$
and 
$$
V_{2,q}=\sum_{j\in D}q_j(1-q_j)(\phi_{2}^{(1)}(p_{j},q_{j}))^2-2\sum_{(i,j)\in  D^2, \ i\neq j}q_iq_j\phi _{2}^{(1)}(p_i,q_i)\phi _{2}^{(1)}(p_j,q_j)
$$  

\begin{theorem}\label{thJ22}
\noindent Under the same assumptions as in Theorem \ref{thJ12}, the following central limit theorems hold.

\begin{itemize}
\item[(a)]One sample : 
 as $n\rightarrow +\infty$,
\begin{equation}
\label{thJ12n1}
\sqrt{n}(J(\widehat{\textbf{p}}_n,\textbf{q})-J(\textbf{p},\textbf{q}))\stackrel{\mathcal{D} }{ \rightsquigarrow} \mathcal{N}\left(0, V_{1,p}
 \right),
\end{equation}

\begin{equation}\sqrt{m}(J(\textbf{p},\widehat{\textbf{q}}_m)-J(\textbf{p},\textbf{q}))\stackrel{\mathcal{D} }{ \rightsquigarrow} \mathcal{N}\left(0, V_{2,q} \right), \label{thJ12n2}
\end{equation}
\item[(b)] Two samples : for $(n,m)\rightarrow (+\infty,+\infty)$ and $nm/(n+m)\rightarrow \gamma\in (0,1)$, 
\begin{equation}
\left(\frac{nm}{mV_{1,p}+
	 nV_{2,q}} \right)^{1/2} \left( J(\widehat{\textbf{p}}_n,\widehat{\textbf{q}}_m)-J(\textbf{p},\textbf{q})\right)\stackrel{\mathcal{D} }{ \rightsquigarrow} \mathcal{N}\left(0,1\right)\label{thJ22n1}
\end{equation}
\end{itemize}
\end{theorem}

\subsection{Direct extensions}$ $ \\

\noindent Quite a few number of divergence measures are not symmetrical. Among these non-symmetrical measures are some of the most interesting ones. For  such measures, estimators of the form $J(\widehat{\textbf{p}}_n,\textbf{q})$, $J(\textbf{p},\widehat{\textbf{q}}_m)$ and $J(\widehat{ \textbf{p}}_n,\widehat{\textbf{q}}_m)$ are not equal to $J(\textbf{q},\widehat{\textbf{p}}_n)$, $J(\widehat{\textbf{q}}_m,\textbf{p})$ and $J(\widehat{\textbf{q}}_m,\widehat{\textbf{p}}_n)$ respectively.\\

\noindent In one-sided tests, we have to decide whether the hypothesis $\textbf{p} =\textbf{q}$, for $\textbf{q}$ known and fixed, is true based on data from $\textbf{p}$. In such a case, we may use the statistics one of the statistics $(J(\widehat{\textbf{p}}_n,\textbf{q})$ and  $J(\textbf{q},\widehat{\textbf{p}}_n)$ to perform the tests. We may have information that allows us to prefer one of them. If not, it is better to use both of them, upon the finiteness of both $J(\textbf{p},\textbf{q})$ and $J(\textbf{q},\textbf{p})$, in a symmetrized form as

\begin{equation}
J^{(s)}(\textbf{p},\textbf{q})=\frac{J(\textbf{p},\textbf{q})+J(\textbf{q},\textbf{p})}{2}.
\end{equation}

\noindent The same situation applies when we face double-side tests, i.e., testing $\textbf{p}=\textbf{q}$ from data generated by $\textbf{p}$ et $\textbf{q}$.\\

\noindent \textbf{Asymptotic a.e. efficiency}. 

\begin{theorem}\label{thJs22}Under the same assumptions as in Theorem \ref{thJ12}, the following hold
\begin{itemize}
\item[(a)]One sample :
\begin{equation}
\label{thJs12c1} \limsup_{n\rightarrow +\infty} \frac{\left\vert J^{(s)}(\widehat{\textbf{p}}_n,\textbf{q})-J^{(s)}(\textbf{p},\textbf{q})\right \vert}{a_{n}} \leq \frac{1}{2}\left( A_{1,p}+A_{4,p}\right)\ \  \text{ a.e.},
\end{equation}
\begin{equation}
\limsup_{n\rightarrow +\infty} \frac{\left\vert J^{(s)}(\textbf{p},\widehat{\textbf{q}}_m)-J^{(s)}(\textbf{p},\textbf{q})\right \vert}{b_{n}} \leq \frac{1}{2}\left( A_{2,q}+A_{3,q} \right)\ \  \text{ a.e.},
\label{thJs22c1}
\end{equation}
\item[(b)]Two samples :
\begin{equation}  \limsup_{(n,m)\rightarrow (+\infty,+\infty)} \frac{\left \vert J^{(s)}(\widehat{\textbf{p}}_n,\widehat{\textbf{q}}_m)-J^{(s)}(\textbf{p},\textbf{q})\right \vert}{c_{n,m}}\leq \frac{1}{2}\left( A_{1,p}+A_{2,q}+A_{3,q}+A_{4,p}\right), \text{ a.e.}
\end{equation}

\end{itemize}\end{theorem}

\bigskip
\noindent \textbf{Asymptotic Normality}. \\

\noindent Denote

\begin{eqnarray*}
V_{3,q}&=&\sum_{j\in D}q_j(1-q_j)(\phi_{1}^{(1)}(q_{j},p_{j}))^2-2\sum_{(i,j)\in  D^2, \ i\neq j}q_iq_j\phi_{1}^{(1)}(q_{i},p_{i})\phi_{1}^{(1)}(q_{j},p_{j}),\\
V_{4,p}&=&\sum_{j\in D}p_j(1-p_j)(\phi _{2}^{(1)}(q_j,p_j))^2-2\sum_{(i,j)\in  D^2, \ i\neq j}p_ip_j\phi_{2}^{(1)}(q_{i},p_{i})\phi_{2}^{(1)}(q_{j},p_{j}),
\end{eqnarray*}
\noindent and finally
$$
V_{1,4,p}=\frac{1}{4}(V_{1,p}+V_{4,p}) \ \text{ and}\  \ V_{2,3,q}=\frac{1}{4}(V_{2,q}+V_{3,q}).
$$

\noindent We have
\begin{theorem} \label{thJs22d} Under the same assumptions as in Theorem \ref{thJ12}, the following hold.
\begin{itemize}
\item[(a)]One sample : as $n\rightarrow+\infty$,
\begin{equation}
\sqrt{\frac{n}{V_{1,4,p}}}  \biggr(J^{(s)}(\widehat{\textbf{p}}_n,\textbf{q})-J^{(s)}(\textbf{p},\textbf{q})\biggr) \stackrel{\mathcal{D} }{ \rightsquigarrow} \mathcal{N}(0, 1),
\end{equation}

\begin{equation}
\sqrt{\frac{n}{V_{2,3,q}}}  \biggr(J^{(s)}(\textbf{p},\widehat{\textbf{q}}_n)-J^{(s)}(\textbf{p},\textbf{q})\biggr) \stackrel{\mathcal{D} }{ \rightsquigarrow} \mathcal{N}(0, 1).
\end{equation}
\item[(b)]Two samples : for $(n,m)\rightarrow (+\infty,+\infty)$ and $nm/(n+m)\rightarrow \gamma\in (0,1)$, 

\begin{equation}
\left(\frac{nm}{mV_{1,4,p}+nV_{2,3,q}}\right)^{1/2} \biggr(J^{(s)}(\widehat{\textbf{p}}_n,\widehat{\textbf{q}}_m)-J^{(s)}(\textbf{p},\textbf{q})\biggr) \stackrel{\mathcal{D} }{ \rightsquigarrow} \mathcal{N}(0, 1).
\end{equation}
\end{itemize}\end{theorem}

\bigskip \noindent \textbf{Remark} : The proof of these extensions will not be given here, since they are straight consequences of the main results. As well, such considerations will not be made again for particular measures for the same reason.
\section{ Particular Cases}\label{particular-divergencemeasure}

\subsection{Renyi and Tsallis families}$ $ \\

\noindent These two families are expressed through the summation
\begin{equation}\label{ialpha}\mathcal{S}_\alpha(\textbf{p},\textbf{q})=\sum_{j\in D}p_j^\alpha q_j^{1-\alpha},\ \ \alpha>0,\ \ \alpha\neq 1,\end{equation}which is of the form of the $\phi-$divergence measure with 
$$\phi(x,y)=x^\alpha y^{1-\alpha},\ \ (x,y)\in \{(p_j,q_j),\ j\in D\}.$$

\bigskip \noindent \textbf{A-(a)- The asymptotic behavior of the Tsallis divergence measure.}$ $ \\ 

\noindent Denote \begin{eqnarray*}
A_{T,\alpha,1}=:\frac{\alpha}{|\alpha-1|}\sum_{j\in D}\left( p_j/q_j\right)^{\alpha-1}\ \ \text{and }\ 
A_{T,\alpha,2}=:\sum_{j\in D}\left( p_j/q_j\right)^{\alpha}.
\end{eqnarray*}
We have 
\begin{corollary}Under the same assumptions as in Theorem \ref{thJ12}, and for any $\alpha>0,\ \alpha\neq 1$, the following hold 
\begin{itemize}
\item[(a)]One sample :
\begin{eqnarray*}
&& \limsup_{n\rightarrow +\infty}\frac{|\mathcal{D}_{T,\alpha}(\widehat{\textbf{p}}_n,\textbf{q})-\mathcal{D}_{T,\alpha}(\textbf{p},\textbf{q})|}{a_n}\leq A_{T,\alpha,1}\ \ \text{a.s} ,\\
&& \limsup_{n\rightarrow +\infty}\frac{|\mathcal{D}_{T,\alpha}(\textbf{p},\widehat{\textbf{q}}_n)-\mathcal{D}_{T,\alpha}(\textbf{p},\textbf{q})|}{b_n}\leq A_{T,\alpha,2}\ \ \text{a.s}.
\end{eqnarray*}
\item[(b)] Two samples :
\begin{eqnarray*}
 \limsup_{(n,m)\rightarrow (+\infty, +\infty)}\frac{|\mathcal{D}_{T,\alpha}(\widehat{\textbf{p}}_n,\widehat{\textbf{q}}_m)-\mathcal{D}_{T,\alpha}(\textbf{p},\textbf{q})|}{c_{n,m}}\leq A_{T,\alpha,1}+A_{T,\alpha,2}\ \ \text{a.s}.
\end{eqnarray*}
\end{itemize}
\end{corollary}

\bigskip \noindent Denote 
 \begin{eqnarray*}
 V_{T,\alpha,1}(\textbf{p},\textbf{q})&=&\left(\frac{\alpha}{\alpha-1}\right)^2\biggr(  \sum_{j\in D}p_j(1-p_j)( p_j/q_j)^{2\alpha-2}\\
 &&\ \ \ \ \ \ \  \ \ \  -\ \ \ 2 \sum_{(i,j)\in D^2,\,i\neq j}(p_ip_j)^\alpha (q_iq_j)^{\alpha-1}\biggr)
 \\
\text{and}\ \  V_{T,\alpha,2}(\textbf{p},\textbf{q})&=& \sum_{j\in D}q_j(1-q_j)(p_j/q_j)^{2\alpha}-2\sum_{(i,j)\in  D^2, \ i\neq j}(p_ip_j)^\alpha (q_iq_j)^{1-\alpha} .
 \end{eqnarray*}

\noindent We have 
\begin{corollary}Under the same assumptions as in Theorem \ref{thJ12}, and for any $\alpha>0,\ \alpha\neq 1$, the following hold 
\begin{itemize}
\item[(a)] One sample : as $n\rightarrow +\infty$
\begin{eqnarray*}
\sqrt{n}\left(\mathcal{D}_{T,\alpha}(\widehat{\textbf{p}}_n,\textbf{q})-\mathcal{D}_{T,\alpha}(\textbf{p},\textbf{q}) \right)\stackrel{\mathcal{D} }{ \rightsquigarrow} \mathcal{N}(0,V_{T,\alpha,1}(\textbf{p},\textbf{q})),
\end{eqnarray*}
\begin{eqnarray*}
\sqrt{n}\left(\mathcal{D}_{T,\alpha}(\textbf{p},\widehat{\textbf{q}}_n)-\mathcal{D}_{T,\alpha}(\textbf{p},\textbf{q}) \right)\stackrel{\mathcal{D} }{ \rightsquigarrow} \mathcal{N}(0,V_{T,\alpha,2}(\textbf{p},\textbf{q})).
\end{eqnarray*}
\item[(b)] Two samples : for $(n,m)\rightarrow (+\infty,+\infty)$ and $nm/(n+m)\rightarrow \gamma\in (0,1)$, 
\begin{eqnarray*}
\left( \frac{mn}{nV_{T,\alpha,2}(\textbf{p},\textbf{q})+mV_{T,\alpha,1}(\textbf{p},\textbf{q}) }\right)^{1/2}\left(\mathcal{D}_{T,\alpha}(\widehat{\textbf{p}}_n,\widehat{\textbf{q}}_m)-\mathcal{D}_{T,\alpha}(\textbf{p},\textbf{q}) \right)\stackrel{\mathcal{D} }{ \rightsquigarrow} \mathcal{N}(0,1).
\end{eqnarray*}
\end{itemize}

\end{corollary}

\bigskip \noindent As to the symmetrized form
 $$ \mathcal{D}_{T,\alpha}^{(s)}(\textbf{p},\textbf{q})=\frac{\mathcal{D}_{T,\alpha}(\textbf{p},\textbf{q})+\mathcal{D}_{T,\alpha}(\textbf{q},\textbf{p})}{2},$$
\noindent we need the supplementaries notations:
\begin{eqnarray*}
A_{T,\alpha,3}&=&\frac{ \alpha}{|\alpha-1|}\sum_{j\in D}\left( q_j/p_j\right)^{\alpha-1},\ \ \ A_{T,\alpha,4}=\sum_{j\in D}\left( q_j/p_j\right)^{\alpha},
\end{eqnarray*}
\begin{eqnarray*}
 V_{T,\alpha,3}(\textbf{p},\textbf{q})&=& \left(\frac{\alpha}{\alpha-1}\right)^2\left( \sum_{j\in D} q_j(1-q_j)(q_j/p_j)^{2-2\alpha}-2\sum_{(i,j)\in  D^2, \ i\neq j}(q_iq_j)^{2-\alpha} (p_ip_j)^{\alpha-1}\right),\\
V_{T,\alpha,4}(\textbf{p},\textbf{q})&=& \sum_{j\in D}p_j(1-p_j) (q_{j}/p_{j})^{2\alpha}- 2 \sum_{(i,j)\in D^2,\,i\neq j}(p_ip_j)^{1-\alpha}(q_{i}q_j)^{\alpha}.
 \end{eqnarray*}
We have 
\begin{corollary}
Under the same assumptions as in Theorem \ref{thJ12}, and for any $\alpha>0,\ \alpha\neq 1$, the following hold 
\begin{itemize}
\item[(a)]One sample :
\begin{eqnarray*}
 \limsup_{n\rightarrow +\infty}\frac{|\mathcal{D}_{T,\alpha}^{(s)}(\widehat{\textbf{p}}_n,\textbf{q})-\mathcal{D}_{T,\alpha}^{(s)}(\textbf{p},\textbf{q})|}{a_n}\leq \left(A_{T,\alpha,1}+A_{T,\alpha,4}\right)/2=:A_{T,\alpha,1}^{(s)}\ \ \text{a.s} ,
\end{eqnarray*}
\begin{eqnarray*}
\limsup_{n\rightarrow +\infty}\frac{|\mathcal{D}_{T,\alpha}^{(s)}( \textbf{p},\widehat{\textbf{q}}_n)-\mathcal{D}_{T,\alpha}^{(s)}(\textbf{p},\textbf{q})|}{b_n}\leq \left(A_{T,\alpha,2}+A_{T,\alpha,3}\right)/2=:A_{T,\alpha,2}^{(s)}\ \ \text{a.s}.
\end{eqnarray*}
\item[(b)]Two samples : \begin{eqnarray*}
 \limsup_{(n,m)\rightarrow (+\infty,+\infty)}\frac{|\mathcal{D}_{T,\alpha}^{(s)}(\widehat{\textbf{p}}_n,\widehat{\textbf{q}}_m)-\mathcal{D}_{T,\alpha}^{(s)}(\textbf{p},\textbf{q})|}{c_{n,m}}\leq A_{T,\alpha,1}^{(s)}+A_{T,\alpha,2}^{(s)}\ \ \text{a.s}.
\end{eqnarray*}
\end{itemize}

\end{corollary}

\bigskip \noindent Denote 
\begin{eqnarray*}
V_{T,\alpha,1:4}(\textbf{p},\textbf{q})&=&V_{T,\alpha,1}(\textbf{p},\textbf{q})+V_{T,\alpha,4}(\textbf{p},\textbf{q})\\
\text{and}\ \  V_{T,\alpha,2:3}(\textbf{p},\textbf{q})&=&V_{T,\alpha,2}(\textbf{p},\textbf{q})+V_{T,\alpha,3}(\textbf{p},\textbf{q}).
\end{eqnarray*}We also have

\begin{corollary}
Under the same assumptions as in Theorem \ref{thJ12}, and for any $\alpha>0,\ \alpha\neq 1$, the following hold
\begin{itemize}
\item[(a)] One sample : as $n\rightarrow +\infty,$
\begin{eqnarray*}
\sqrt{n}\biggr(\mathcal{D}_{T,\alpha}^{(s)}(\widehat{\textbf{p}}_n,\textbf{q})-\mathcal{D}_{T,\alpha}^{(s)}(\textbf{p},\textbf{q})\biggr) \stackrel{\mathcal{D} }{ \rightsquigarrow} \mathcal{N}(0,V_{T,\alpha,1:4}(\textbf{p},\textbf{q})),
\end{eqnarray*}
\begin{eqnarray*}
\sqrt{n}\biggr(\mathcal{D}_{T,\alpha}^{(s)}(\textbf{p},\widehat{\textbf{q}}_n)-\mathcal{D}_{T,\alpha}^{(s)}(\textbf{p},\textbf{q})\biggr) \stackrel{\mathcal{D} }{ \rightsquigarrow} \mathcal{N}(0,V_{T,\alpha,2:3}(\textbf{p},\textbf{q})).
\end{eqnarray*}
\item[(b)] Two samples :  for $(n,m)\rightarrow (+\infty,+\infty)$ and $nm/(n+m)\rightarrow \gamma\in (0,1)$, 
\begin{eqnarray*}
&&\left(\frac{nm}{mV_{T,\alpha,1:4}(\textbf{p},\textbf{q})+nV_{T,\alpha,2:3}(\textbf{p},\textbf{q})}\right)^{1/2} \biggr(\mathcal{D}_{T,\alpha}^{(s)}(\widehat{\textbf{p}}_n,\widehat{\textbf{q}}_m)-\mathcal{D}_{T,\alpha}^{(s)}(\textbf{p},\textbf{q})\biggr) \stackrel{\mathcal{D} }{ \rightsquigarrow} \mathcal{N}(0, 1).
\end{eqnarray*}
\end{itemize}
\end{corollary}

\bigskip \noindent \textbf{A-(b)- The asymptotic behavior of the Renyi-$\alpha$ divergence measure.}$ $ \\ 

\noindent The treatment of the asymptotic behavior of the Renyi-$\alpha$, $\alpha>0$, $\alpha\neq 1$ is obtained from Part \textbf{(A)-(a)} by expansions and by the application of the delta method.\\

\noindent We first remark that 
\begin{eqnarray*}
\mathcal{D}_{R,\alpha}(\textbf{p},\textbf{q})=\frac{1}{\alpha-1}\log \left( \mathcal{S}_\alpha(\textbf{p},\textbf{q})\right).
\end{eqnarray*}
\begin{corollary}
Under the same assumptions as in Theorem \ref{thJ12}, and for any $\alpha>0,\ \alpha\neq 1$, the following hold
\begin{itemize}
\item[(a)]One sample :
\begin{eqnarray*}
\limsup_{n\rightarrow +\infty}\frac{|\mathcal{D}_{R,\alpha}(\widehat{\textbf{p}}_n,\textbf{q})-\mathcal{D}_{R,\alpha}(\textbf{p},\textbf{q})|}{a_n}\leq \frac{ A_{T,\alpha,1}}{\mathcal{S}_\alpha(\textbf{p},\textbf{q})}=:A_{R,\alpha,1}\ \ \text{a.s},
\end{eqnarray*}
\begin{eqnarray*}
\limsup_{n\rightarrow +\infty}\frac{|\mathcal{D}_{R,\alpha}( \textbf{p},\widehat{\textbf{q}}_n)-\mathcal{D}_{R,\alpha}(\textbf{p},\textbf{q})|}{b_n}\leq\frac{ A_{T,\alpha,2}}{\mathcal{S}_\alpha(\textbf{p},\textbf{q})}=:A_{R,\alpha,2}\ \ \text{a.s}.
\end{eqnarray*}
\item[(b)]Two samples :
\begin{eqnarray*}
  \limsup_{(n,m)\rightarrow( +\infty,+\infty)}\frac{|\mathcal{D}_{R,\alpha}(\widehat{\textbf{p}}_n,\widehat{\textbf{q}}_m)-\mathcal{D}_{R,\alpha}(\textbf{p},\textbf{q})|}{c_{n,m}}\leq A_{R,\alpha,1}+A_{R,\alpha,2}\ \ \text{a.s}.
\end{eqnarray*}
\end{itemize}
\end{corollary}

\bigskip \noindent Denote 
$$V_{R,\alpha,1}(\textbf{p},\textbf{q})=\frac{V_{T,\alpha,1}(\textbf{p},\textbf{q})}{\mathcal{S}_\alpha^2(\textbf{p},\textbf{q})}\ \ \text{and}\ \ V_{R,\alpha,2}(\textbf{p},\textbf{q})=\frac{V_{T,\alpha,2}(\textbf{p},\textbf{q})}{\mathcal{S}_\alpha^2(\textbf{p},\textbf{q})}.$$ We have \begin{corollary}Under the same assumptions as in Theorem \ref{thJ12}, and for any $\alpha>0,\ \alpha\neq 1$, the following hold
\begin{itemize}
\item[(a)]One sample : as $n\rightarrow+\infty$
\begin{eqnarray*}
\sqrt{n}\biggr(\mathcal{D}_{R,\alpha}(\widehat{\textbf{p}}_n,\textbf{q})-\mathcal{D}_{R,\alpha}(\textbf{p},\textbf{q})\biggr)\stackrel{\mathcal{D} }{ \rightsquigarrow} \mathcal{N}(0,V_{R,\alpha,1}(\textbf{p},\textbf{q})),
\end{eqnarray*}
\begin{eqnarray*}
\sqrt{n}\biggr(\mathcal{D}_{R,\alpha}(\textbf{p},\widehat{\textbf{q}}_n)-\mathcal{D}_{R,\alpha}(\textbf{p},\textbf{q})\biggr)\stackrel{\mathcal{D} }{ \rightsquigarrow} \mathcal{N}(0,V_{R,\alpha,2}(\textbf{p},\textbf{q})).
\end{eqnarray*}
\item[(b)]Two samples : for $(n,m)\rightarrow (+\infty,+\infty)$ and $nm/(n+m)\rightarrow \gamma\in (0,1)$, 
\begin{eqnarray*}
\left(\frac{mn}{nV_{R,\alpha,2}(\textbf{p},\textbf{q})+mV_{R,\alpha,1}(\textbf{p},\textbf{q})}\right)^{1/2}\biggr(\mathcal{D}_{R,\alpha}(\widehat{\textbf{p}}_n,\widehat{\textbf{q}}_m)-\mathcal{D}_{R,\alpha}(\textbf{p},\textbf{q})\biggr)\stackrel{\mathcal{D} }{ \rightsquigarrow} \mathcal{N}(0,1).
\end{eqnarray*}\end{itemize}

\end{corollary}

\bigskip \noindent As to the symetrized form
$$\mathcal{D}_{R,\alpha}^{(s)}(\textbf{p},\textbf{q})=\frac{\mathcal{D}_{R,\alpha}(\textbf{p},\textbf{q})-\mathcal{D}_{R,\alpha}(\textbf{q},\textbf{p})}{2},$$

\noindent we need the supplementary notations :

\begin{eqnarray*}
A_{R,\alpha,3}&=& \frac{ A_{T,\alpha,3}}{\mathcal{S}_\alpha(\textbf{p},\textbf{q})},\ \ \ \ A_{R,\alpha,4}=\frac{ A_{T,\alpha,4}}{\mathcal{S}_\alpha(\textbf{p},\textbf{q})}\\
V_{R,\alpha,3}(\textbf{p},\textbf{q})&=&\frac{V_{T,\alpha,3}(\textbf{p},\textbf{q})}{\mathcal{S}_\alpha^2(\textbf{p},\textbf{q})}\ \ \text{and}\ \ V_{R,\alpha,4}(\textbf{p},\textbf{q})=\frac{V_{T,\alpha,4}(\textbf{p},\textbf{q})}{\mathcal{S}_\alpha^2(\textbf{p},\textbf{q})}.
\end{eqnarray*}
\begin{corollary}Under the same assumptions as in Theorem \ref{thJ12}, and for any $\alpha>0,\ \alpha\neq 1$, the following hold.
\begin{itemize}
\item[(a)] One sample :
\begin{eqnarray*}
&& \limsup_{n\rightarrow+\infty}\frac{|\mathcal{D}_{R,\alpha}^{(s)}(\widehat{\textbf{p}}_n,\textbf{q})-\mathcal{D}_{R,\alpha}^{(s)}(\textbf{p},\textbf{q})|}{a_n}\leq (A_{R,\alpha,1}+A_{R,\alpha,4})/2=:A_{R,\alpha,1}^{(s)},\ \ \text{a.s.}\\
&& \limsup_{n\rightarrow+\infty}\frac{|\mathcal{D}_{R,\alpha}^{(s)}(\textbf{p},\widehat{\textbf{q}}_n)-\mathcal{D}_{R,\alpha}^{(s)}(\textbf{p},\textbf{q})|}{a_n}\leq (A_{R,\alpha,2}+A_{R,\alpha,3})/2=:A_{R,\alpha,2}^{(s)}.
\end{eqnarray*}
\item[(b)] Two samples :
\begin{eqnarray*}
 \limsup_{(n,m)\rightarrow(+\infty,+\infty)}\frac{|\mathcal{D}_{R,\alpha}^{(s)}(\widehat{\textbf{p}}_n,\widehat{\textbf{q}}_m)-\mathcal{D}_{R,\alpha}^{(s)}(\textbf{p},\textbf{q})|}{c_{n,m}}\leq A_{R,\alpha,1}^{(s)}+A_{R,\alpha,2}^{(s)},\ \ \text{a.s.}
\end{eqnarray*}
\end{itemize}
\end{corollary}
\noindent Denote 
\begin{eqnarray*}
V_{R,\alpha,1:4}(\textbf{p},\textbf{q})&=&V_{R,\alpha,1}(\textbf{p},\textbf{q})+V_{R,\alpha,4}(\textbf{p},\textbf{q})\\
\text{and}\ \ V_{R,\alpha,2:3}(\textbf{p},\textbf{q})&=&V_{R,\alpha,2}(\textbf{p},\textbf{q})+V_{R,\alpha,3}(\textbf{p},\textbf{q}),\ \ \text{a.s.}
\end{eqnarray*}

\bigskip \noindent We also have 
\begin{corollary}
Under the same assumptions as in Theorem \ref{thJ12}, and for any $\alpha>0,\ \alpha\neq 1$, the following hold.
\begin{itemize}
\item[(a)]One sample : as $n\rightarrow+\infty$,
\begin{eqnarray*}
&&\sqrt{n}\biggr(\mathcal{D}_{R,\alpha}^{(s)}(\widehat{\textbf{p}}_n,\textbf{q})-\mathcal{D}_{R,\alpha}^{(s)}(\textbf{p},\textbf{q}) \biggr)\stackrel{\mathcal{D} }{ \rightsquigarrow} \mathcal{N}(0,V_{R,\alpha,1:4}(\textbf{p},\textbf{q})),\ \ \text{as} \ \ n\rightarrow +\infty\\
&&\sqrt{n}\biggr(\mathcal{D}_{R,\alpha}^{(s)}(\textbf{p},\widehat{\textbf{q}}_n)-\mathcal{D}_{R,\alpha}^{(s)}(\textbf{p},\textbf{q}) \biggr)\stackrel{\mathcal{D} }{ \rightsquigarrow} \mathcal{N}(0,V_{R,\alpha,2:3}(\textbf{p},\textbf{q})),\ \ \text{as} \ \ n\rightarrow +\infty
\end{eqnarray*}
\item[(b)] Two samples : as $(n,m)\rightarrow (+\infty,+\infty)$ and $nm/(n+m)\rightarrow \gamma\in (0,1)$, 
\begin{eqnarray*}
&& \left(\frac{mn}{nV_{R,\alpha,2:3}(\textbf{p},\textbf{q})+mV_{R,\alpha,1:4}(\textbf{p},\textbf{q})}\right)^{1/2}\biggr(\mathcal{D}_{R,\alpha}^{(s)}(\widehat{\textbf{p}}_n,\widehat{\textbf{q}}_m)-\mathcal{D}_{R,\alpha}^{(s)}(\textbf{p},\textbf{q}) \biggr)\stackrel{\mathcal{D} }{ \rightsquigarrow} \mathcal{N}(0,1).
\end{eqnarray*}
\end{itemize}
\end{corollary}

\bigskip \textbf{B - Kulback-Leibler Measure}$ $\\

\bigskip \noindent Here we have 
$$
\mathcal{D}_{KL}(\textbf{p},\textbf{q})=\sum_{j\in D}\phi(p_j,q_j),$$ where
$$\phi(x,y)=x\log (x/y),\ \ (x,y)\in \{(p_j,q_j),\,j\in D\}.$$

\noindent So we have first :

\begin{corollary}Under the same assumptions as in Theorem \ref{thJ12}, the following hold.

\begin{itemize}
\item[(a)] One sample :
\begin{eqnarray*}
&&\limsup_{n\rightarrow +\infty}\frac{|\mathcal{D}_{KL}(\widehat{\textbf{p}}_n,\textbf{q})-\mathcal{D}_{KL}(\textbf{p},\textbf{q})|}{a_n}\leq \sum_{j\in D}\left\vert 1+\log (p_j/q_j)\right\vert   =:A_{KL,1}(\textbf{p},\textbf{q}),\ \ \text{a.s.}\\
&&\limsup_{n\rightarrow +\infty}\frac{|\mathcal{D}_{KL}(\textbf{p},\widehat{\textbf{q}}_n)-\mathcal{D}_{KL}(\textbf{p},\textbf{q})|}{b_n}\leq \sum_{j\in D}\left(p_j/q_j\right)  =:A_{KL,2}(\textbf{p},\textbf{q}),\ \ \text{a.s.}.
\end{eqnarray*}
\item[(b)] Two samples :
\begin{eqnarray*}
\limsup_{(n,m)\rightarrow (+\infty,+\infty)}\frac{|\mathcal{D}_{KL}(\widehat{\textbf{p}}_n,\widehat{\textbf{q}}_m)-\mathcal{D}_{KL}(\textbf{p},\textbf{q})|}{c_{n,m}}\leq A_{KL,1}(\textbf{p},\textbf{q})+A_{KL,2}(\textbf{p},\textbf{q}),\ \ \text{a.s.}
\end{eqnarray*}
\end{itemize}
\end{corollary}
Denote 
\begin{eqnarray*}
V_{KL,1}(\textbf{p},\textbf{q})&=&
 \sum_{j\in D}p_j(1-p_j)\left(1+\log(p_j/q_j)\right)^2\\
 &&\ \ \ \ -\ \ 2\sum_{(i,j)\in  D^2, \ i\neq j}p_ip_j( 1+\log(p_i/q_{i}))( 1+\log(p_j/q_{j}))\\
\text{and}\ \ \ V_{KL,2}(\textbf{p},\textbf{q})&=&\sum_{j\in D}q_j(1-q_j)(p_j/q_j)^2-2\sum_{(i,j)\in  D^2, \ i\neq j}p_ip_j.
\end{eqnarray*}
 We have 
\begin{corollary}\label{corkbnorm}
Under the same assumptions as in Theorem \ref{thJ12}, the following hold.
\begin{itemize}
\item[(a)] One sample : as $n \rightarrow +\infty$
\begin{eqnarray*}
&& \sqrt{n}\left(\mathcal{D}_{KL}(\widehat{\textbf{p}}_n,\textbf{q})-\mathcal{D}_{KL}(\textbf{p},\textbf{q}) \right)\stackrel{\mathcal{D} }{ \rightsquigarrow} \mathcal{N}(0,V_{KL,1}(\textbf{p},\textbf{q})),\\
&& \sqrt{n}\left(\mathcal{D}_{KL}(\textbf{p},\widehat{\textbf{q}}_n)-\mathcal{D}_{KL}(\textbf{p},\textbf{q}) \right)\stackrel{\mathcal{D} }{ \rightsquigarrow} \mathcal{N}(0,V_{KL,2}(\textbf{p},\textbf{q})).
\end{eqnarray*}
\item[(b)] Two samples : as $(n,m)\rightarrow (+\infty,+\infty)$ and $nm/(n+m)\rightarrow \gamma\in (0,1)$, 
\begin{eqnarray*}
 \left(\frac{mn}{n V_{KL,2}(\textbf{p},\textbf{q})+mV_{KL,1}(\textbf{p},\textbf{q}) } \right)^{1/2}\left(\mathcal{D}_{KL}(\widehat{\textbf{p}}_n,\widehat{\textbf{q}}_m)-\mathcal{D}_{KL}(\textbf{p},\textbf{q}) \right)\stackrel{\mathcal{D} }{ \rightsquigarrow} \mathcal{N}(0,1).
\end{eqnarray*}
\end{itemize}
\end{corollary}

\bigskip \noindent As to the symmetrized form 

$$\mathcal{D}_{KL}^{(s)}(\textbf{p},\textbf{q})=\frac{\mathcal{D}_{KL}(\textbf{p},\textbf{q})+\mathcal{D}_{KL}(\textbf{q},\textbf{p})}{2},$$ 

 \noindent we need the supplementary notations :

\begin{eqnarray*}
 A_{KL,3}(\textbf{p},\textbf{q})=\sum_{j\in D}\left\vert 1+\log(q_j/p_j)\right \vert ,\ \ \ A_{KL,4}(\textbf{p},\textbf{q})=\sum_{j\in D}q_j/p_j.
\end{eqnarray*}

\noindent We have

\begin{corollary}
Under the same assumptions as in Theorem \ref{thJ12}, the following hold
\begin{itemize}
\item[(a)] One sample :
\begin{eqnarray*}
&&\limsup_{n\rightarrow +\infty}\frac{|\mathcal{D}_{KL}^{(s)}(\widehat{\textbf{p}}_n,\textbf{q})-\mathcal{D}_{KL}^{(s)}(\textbf{p},\textbf{q})|}{a_n}\leq  (A_{KL,1}(\textbf{p},\textbf{q})+A_{KL,4}(\textbf{p},\textbf{q}))/2 =:A_{KL,1}^{(s)}(\textbf{p},\textbf{q}),\ \ \text{a.s.}\\
&&\limsup_{n\rightarrow +\infty}\frac{|\mathcal{D}_{KL}^{(s)}(\textbf{p},\widehat{\textbf{q}}_n)-\mathcal{D}_{KL}^{(s)}(\textbf{p},\textbf{q})|}{b_n}\leq (A_{KL,2}(\textbf{p},\textbf{q})+A_{KL,3}(\textbf{p},\textbf{q}))/2 =:A_{KL,2}^{(s)}(\textbf{p},\textbf{q}),\ \ \text{a.s.}
\end{eqnarray*}
\item[(b)]Two samples :
\begin{eqnarray*}
\limsup_{(n,m)\rightarrow (+\infty,+\infty)}\frac{|\mathcal{D}_{KL}^{(s)}(\widehat{\textbf{p}}_n,\widehat{\textbf{q}}_m)-\mathcal{D}_{KL}^{(s)}(\textbf{p},\textbf{q})|}{c_{n,m}}\leq A_{KL,1}^{(s)}(\textbf{p},\textbf{q})+A_{KL,2}^{(s)}(\textbf{p},\textbf{q}),\ \ \text{a.s.}
\end{eqnarray*}
\end{itemize}
\end{corollary}
\noindent Denote 
\begin{eqnarray*}
V_{KL,3}(\textbf{p},\textbf{q})&=&\sum_{j\in D}q_j(1-q_j)(1+\log(q_j/p_j))^2\\
&&\ \ \ \ -2\sum_{(i,j)\in  D^2, \ i\neq j}q_iq_j(1+\log(q_i/p_i))(1+\log(q_j/p_j)),\ \\
V_{KL,4}(\textbf{p},\textbf{q})&=&\sum_{j\in D}p_j(1-p_j)(q_j/p_j)^2-2\sum_{(i,j)\in  D^2, \ i\neq j}q_iq_j,
\end{eqnarray*}and finally
\begin{eqnarray*}
V_{KL,1:4}(\textbf{p},\textbf{q})&=&V_{KL,1}(\textbf{p},\textbf{q})+V_{KL,4}(\textbf{p},\textbf{q}),\\
V_{KL,2:3}(\textbf{p},\textbf{q})&=&V_{KL,2}(\textbf{p},\textbf{q})+V_{KL,3}(\textbf{p},\textbf{q}).
\end{eqnarray*}

\bigskip \noindent We also have 
\begin{corollary}
Under the same assumptions as in Theorem \ref{thJ12}, the following hold
\begin{itemize}
\item[(a)]One sample : as $n \rightarrow +\infty$
\begin{eqnarray*}
&&\sqrt{n}\biggr(\mathcal{D}_{KL}^{(s)}(\widehat{\textbf{p}}_n,\textbf{q})-\mathcal{D}_{KL}^{(s)}(\textbf{p},\textbf{q})\biggr) \stackrel{\mathcal{D} }{ \rightsquigarrow} \mathcal{N}(0, V_{KL,1:4}(\textbf{p},\textbf{q})),
\\
&& \sqrt{n}\biggr(\mathcal{D}_{KL}^{(s)}(\textbf{p},\widehat{\textbf{q}}_n)-\mathcal{D}_{KL}^{(s)}(\textbf{p},\textbf{q})\biggr) \stackrel{\mathcal{D} }{ \rightsquigarrow} \mathcal{N}(0, V_{KL,2:3}(\textbf{p},\textbf{q})).
\end{eqnarray*}
\item[(b)]Two samples : for $(n,m)\rightarrow (+\infty,+\infty)$ and $nm/(n+m)\rightarrow \gamma\in (0,1)$, 

\begin{eqnarray*}
&&\left(\frac{nm}{mV_{KL,1:4}(\textbf{p},\textbf{q})+nV_{KL,2:3}(\textbf{p},\textbf{q})}\right)^{1/2} \biggr(\mathcal{D}_{KL}^{(s)}(\widehat{\textbf{p}}_n,\widehat{\textbf{q}}_m)-\mathcal{D}_{KL}^{(s)}(\textbf{p},\textbf{q})\biggr) \stackrel{\mathcal{D} }{ \rightsquigarrow} \mathcal{N}(0, 1).
\end{eqnarray*}
\end{itemize}
\end{corollary}
\section{Proofs}
\label{proofs}
\noindent In the proofs, we will systematically use the mean values theorem.  In the multivariate handling, we prefer to use the Taylor-Lagrange-Cauchy as stated
in \cite{valiron}, page 230.\\

\noindent 
For sake of simplicity, we introduce the two following notations : $$\Delta_{p_n}^{c_j}=\widehat{p}_n^{c_j}-p_j\ \ \text{and}\ \  \Delta_{q_m}^{c_j}=\widehat{q}_m^{c_j}-q_j,\ \ \forall\,j\in D,$$ therefore $$a_n=\sup_{j\in D}|\Delta_{p_n}^{c_j}|,\ \ \ \ b_m=\sup_{j\in D}|\Delta_{q_m}^{c_j}|,\ \ 
\text{and}\ \ c_{n,m}=\max(a_n,b_m).$$
For any $j\in D$, set
$$\delta_n(p_j)=\sqrt{n/p_j}\Delta_{p_n}^{c_j}\ \ \text{and}\ \ \delta_m(q_j)=\sqrt{m/q_j}\Delta_{q_m}^{c_j}.$$

\bigskip \noindent Before we start the proofs we recall that,
since for a fixed $j\in D,$ $n\widehat{p}_n^{c_j}$ has a binomial distribution with parameters $n$  and success probability $p_j$, we have 
 \begin{equation*}
 \mathbb{E}\left[ \widehat{p}_n^{c_j}\right]=p_j\ \ \text{and}\ \ \mathbb{V}(\widehat{p}_n^{c_j})=\frac{p_j(1-p_j)}{n}.
\end{equation*} Furthermore, by the strong law of large numbers, we know that
$$\Delta_{p_n}^{c_j}\stackrel{a.s.}{\longrightarrow} 0,\ \ \text{as }\ \  n  \rightarrow+\infty,$$
for a fixed $j\in D.$\\

\noindent And finally, by the asymptotic Gaussian limit of the multinomial law (see for example \cite{ips-wcia-ang}, Chapter 1, Section 4), we have
\begin{eqnarray}
\label{pnj}&& \biggr( \delta_n(p_j), \ j\in D\biggr)
\stackrel{\mathcal{D}}{\rightsquigarrow }Z(\textbf{p}) \stackrel{\mathcal{L} }{\sim}\mathcal{N}(0,\Sigma_\textbf{p}),\ \ \text{as}\ \ n\rightarrow +\infty,\\
 \text{and}\ \  && \biggr( \delta_m(q_j), \ j\in D)\biggr)\stackrel{\mathcal{D}}{\rightsquigarrow }Z(\textbf{q})\stackrel{\mathcal{L} }{\sim}\mathcal{N}(0,\Sigma_\textbf{q}),\ \ \text{as}\ \ m\rightarrow +\infty,
 \label{qnj}
\end{eqnarray}where $Z(\textbf{p})= (Z_{p_j},j\in D)^t
$ \
and $Z(\textbf{q})=(Z_{q_j},j\in D)^t
$ are two centered Gaussian random vectors of dimension $\#(D)$
  which are independent and have the following elements :
\begin{eqnarray}\label{vars}
&&\left(\Sigma_\textbf{p}\right)_{(i,j)}=(1-p_j)1_{(i=j)}-\sqrt{ p_ip_j} 1_{(i\neq j)}, \ \ (i,j) \in D^2\\
&&\left(\Sigma_\textbf{q}\right)_{(i,j)}=(1-q_j)1_{(i=j)}-\sqrt{q_iq_j} 1_{(i\neq j)}, \ \ (i,j) \in D^2.
\end{eqnarray}

\bigskip \noindent For a fixed $j\in D$, we have also 
\begin{eqnarray*}
&&\mathbb{E}\left[ \widehat{q}_m^{c_j}\right]=q_j,\ \ \ \ \mathbb{V}(\widehat{q}_m^{c_j})=\frac{q_j(1-q_j)}{m},
\ \ \ \ \Delta_{q_m}^{c_j}\stackrel{a.s.}{\longrightarrow} 0,\ \text{as}\ m\rightarrow+\infty.
\end{eqnarray*}

\bigskip
\noindent Now we can start by showing Theorem \ref{thJ12}. 

\bigskip

\noindent For a fixed $j\in D$, we have 
\begin{eqnarray}\label{fi1}
\phi(\widehat{ p}_n^{c_j},q_j)&=&\phi(p_j+\Delta_{p_n}^{c_j},q_j)\\
&=&\phi(p_j,q_j)+\Delta_{p_n}^{c_j}\phi_{1}^{(1)}(p_{j}+\theta _{1,j}\Delta_{p_n}^{c_j},q_{j})\notag
\end{eqnarray}

\noindent by applying the mean value theorem to the function $(.)\mapsto
\phi((.),q_{j})$ and where $\theta _{1,j}$ is some number lying between $0$ and $1$. In the sequel, any $ \theta_{i,j},\,i=1,2,\cdots$ satisfies $ \left| \theta_{i,j} \right|<1$. \\

\noindent By applying again the mean values theorem to the function $(.)\mapsto
\phi_{1}^{(1)}((.),q_{j})$, we have
\begin{eqnarray*}\phi_{1}^{(1)}(p_{j}+\theta _{1,j}\Delta_{p_n}^{c_j},q_{j})&=&\phi_{1}^{(1)}(p_{j},q_{j})+\theta _{1,j}\Delta_{p_n}^{c_j}\phi_{1}^{(2)}(p_{j}+\theta _{2,j}\Delta_{p_n}^{c_j},q_{j})
\end{eqnarray*}
\noindent We can write \eqref{fi1} as 
\begin{eqnarray*}
\phi(\widehat{ p}_n^{c_j},q_j)&=&\phi(p_j,q_j)+\Delta_{p_n}^{c_j} \phi_{1}^{(1)}(p_{j},q_{j})+\theta _{1,j}(\Delta_{p_n}^{c_j})^2\phi_{1}^{(2)}(p_{j}+\theta _{2,j}\Delta_{p_n}^{c_j},q_{j})
\end{eqnarray*}

\noindent Now we have
\begin{eqnarray}\label{jpn}
\notag J(\widehat{\textbf{p}}_n,\textbf{q})-J(\textbf{p},\textbf{q})&=&\sum_{j\in D} \Delta_{p_n}^{c_j} \phi_{1}^{(1)}(p_{j},q_{j})\\
&& \ \ +\ \ \sum_{j\in D}\theta _{1,j}(\Delta_{p_n}^{c_j})^2\phi_{1}^{(2)}(p_{j}+\theta _{2,j}\Delta_{p_n}^{c_j},q_{j}),
\end{eqnarray}

\noindent hence \begin{eqnarray*}
\left\vert J(\widehat{\textbf{p}}_n,\textbf{q})-J(\textbf{p},\textbf{q})\right\vert &\leq & a_n\sum_{j\in D} | \phi_{1}^{(1)}(p_{j},q_{j})| \ + \ a_n^2\sum_{j\in D}|\phi_{1}^{(2)}(p_{j}+\theta _{2,j}\Delta_{p_n}^{c_j},q_{j})|,
\end{eqnarray*}
 
\noindent Therefore 
\begin{equation*}
\limsup_{n\rightarrow \infty }\frac{|J(\widehat{\textbf{p}}_n,\textbf{q})-J(\textbf{p},\textbf{q})}{a_{n}}\leq A_{1,p}+
a_{n} \sum_{j\in D}|\phi_{1}^{(2)}(p_{j}+\theta _{2,j}\Delta_{p_n}^{c_j},q_{j})|.
\end{equation*}
We know that $A_{1,p}<\infty$ and 
\begin{equation*}
 \sum_{j\in D}|\phi_{1}^{(2)}(p_{j}+\theta _{2,j}\Delta_{p_n}^{c_j},q_{j})[ \rightarrow  \sum_{j\in D}|\phi_{1}^{(2)}(p_{j},q_{j})|<\infty \ \text{\ \ as \
\ } n\rightarrow \infty.
\end{equation*}
\bigskip
\noindent This proves \eqref{thJ12c1}.\\

\bigskip \noindent Formula \eqref{thJ12c2} is obtained in a similar way. We only need to adapt the result concerning the first coordinate to the second. 
 
 \bigskip \noindent The proof of \eqref{thJ22c1} comes by splitting $\sum_{j\in D}\left(\phi(\widehat{p}_n^{c_j},\widehat{q}_m^{c_j})-\phi(p_j,q_j) \right)$, into the following two terms 

\begin{eqnarray*}
\sum_{j\in D}\left(\phi(\widehat{p}_n^{c_j},\widehat{q}_m^{c_j})-\phi(p_j,q_j) \right)&=& \sum_{j\in D}\left(\phi(\widehat{p}_n^{c_j},\widehat{q}_m^{c_j})-\phi(p_j,\widehat{q}_m^{c_j}) \right)\\
&+& \sum_{j\in D}\left(\phi(p_j,\widehat{q}_m^{c_j})- \phi(p_j,q_j)\right)\\
&\equiv & I_{n,1}+ I_{n,2} 
\end{eqnarray*}

\bigskip \noindent We already know how the handle $I_{n,2}$. As to $I_{n,1}$, we may still use the Taylor-Lagrange-Cauchy formula since we have, for a fixed $j\in D$,

$$
\left\Vert (\widehat{p}_n^{c_j},\widehat{q}_m^{c_j})-(p_j,\widehat{q}_m^{c_j})\right\Vert_{\infty}=\left\Vert (\widehat{p}_n^{c_j}-p_j),0)\right\Vert_{\infty}=a_n\rightarrow 0.
$$

\bigskip \noindent By the Taylor-Lagrange-Cauchy (see \cite{valiron}, page 230), we have

\begin{eqnarray*}
I_{n,1}&=&\sum_{j\in D} \Delta_{p_n}^{c_j} \phi(\widehat{p}_n^{c_j}+\theta_j \Delta_{p_n}^{c_j},\widehat{q}_m^{c_j}) \\
&\leq & a_n \sum_{j\in D}| \phi(\widehat{p}_n^{c_j}+\theta_j \Delta_{p_n}^{c_j},\widehat{q}_m^{c_j})| \\
&=& a_n (A_1 + o(1)).
\end{eqnarray*}

\noindent From there, the combination of these remarks direct to the result.$\blacksquare$\\

\bigskip

\noindent Let us prove \eqref{thJ12n1}. By going back to \eqref{jpn}, we  have 

\begin{eqnarray*}
\sqrt{n}(J(\widehat{\textbf{p}}_n,\textbf{q})-J(\textbf{p},\textbf{q}))&=& \sum_{j\in D}\sqrt{p_j}\delta_n(p_j) \phi_{1}^{(1)}(p_{j},q_{j}) + \sqrt{n}R_{1,n}
 \end{eqnarray*}

\bigskip \noindent where $$ R_{1,n}=\sum_{j\in D}\theta _{1,j}(\Delta_{p_n}^{c_j})^2\phi_{1}^{(2)}(p_{j}+\theta _{2,j}\Delta_{p_n}^{c_j},q_{j}). $$ 

\bigskip \noindent Now using Formula \eqref{pnj} above, we get,
$$
\sum_{j\in D}\sqrt{p_j}\delta_n(p_j) \phi_{1}^{(1)}(p_{j},q_{j})\stackrel{\mathcal{D} }{\rightsquigarrow} \sum_{j\in D} \phi_{1}^{(1)}(p_{j},q_{j})\sqrt{p_j}Z_{p_{j}},\ \ \text{as}\ \ n\rightarrow+\infty
$$
\noindent which follows a centered normal law of variance $V_{1,p}$ : $$V_{1,p}=\sum_{j\in D}(1-p_j)(\phi_{1}^{(1)}(p_{j},q_{j}))^2-2\sum_{(i,j)\in D^2,\,i\neq j}\sqrt{p_ip_j}\phi_{1}^{(1)}(p_{i},q_{i})\phi_{1}^{(1)}(p_{j},q_{j})$$
since 
\begin{eqnarray*}
\text{V}ar\left(\sum_{j\in D}\phi_{1}^{(1)}(p_{j},q_{j})\sqrt{p_j}Z_{p_j}\right)&=&\sum_{j\in D}\text{V}ar(\phi_{1}^{(1)}(p_{j},q_{j})\sqrt{p_j}Z_{p_j})\\
\ \ \ \ \  & & \ \ \ \ +  \ \ 2\sum_{(i,j)\in D^2,\,i\neq j}\text{Cov}(\phi_{1}^{(1)}(p_{i},q_{i})\sqrt{p_i}Z_{p_i},\phi_{1}^{(1)}(p_{j},q_{j})\sqrt{p_j}Z_{p_j})\\
&=&\sum_{j\in D}p_j(\phi_{1}^{(1)}(p_{j},q_{j}))^2\text{V}ar(Z_{p_j})\\
\ \ \ \ \  & & \ \ \ \ +  \ \              2\sum_{(i,j)\in D^2,\,i\neq j}\sqrt{p_ip_j}\phi_{1}^{(1)}(p_{i},q_{i})\phi_{1}^{(1)}(p_{j},q_{j})\text{Cov}(Z_{p_i},Z_{p_j})
\\
&=&\sum_{j\in D}p_j(1-p_j)(\phi_{1}^{(1)}(p_{j},q_{j}))^2\\
\ \ \ \ \  & & \ \ \ \ -  \ \              2\sum_{(i,j)\in D^2,\,i\neq j}p_ip_j\phi_{1}^{(1)}(p_{i},q_{i})\phi_{1}^{(1)}(p_{j},q_{j}).
\end{eqnarray*}

\bigskip \noindent Let show that $\sqrt{n}R_{1,n}=0_{%
\mathbb{P}}(1)$. We have  
\begin{equation} \label{r1n}
\left\vert \sqrt{n}R_{1,n}\right\vert \leq \sqrt{n}a_{n}^{2} \sum_{j\in D}|\phi_{1}^{(2)}(p_{j}+\theta _{2,j}\Delta_{p_n}^{c_j},q_{j})|.
\end{equation}

\bigskip \noindent Let show that $$\sqrt{n}a_{n}^{2}=o_{\mathbb{P}}(1).$$ By the Bienaymé-Tchebychev inequality, we have, for any $\epsilon >0$ and for $j\in D$,
\begin{eqnarray*}
\mathbb{P}(\sqrt{n}(\widehat{p}_n^{c_j}-p_j)^2\geq \epsilon)=\mathbb{P}\left(|\widehat{p}_n^{c_j}-p_j|\geq  \frac{\sqrt{\epsilon} }{n^{1/4}}\right)\leq \frac{p_j(1-p_j)}{\epsilon n^{1/2}},
\end{eqnarray*}which implies that $\sqrt{n}a_{n}^{2}$ converges in probability to $0$ as $n\rightarrow+\infty$.\\

\noindent Finally from \eqref{r1n} we have $\sqrt{n}
R_{1,n}\stackrel{\mathbb{P}}{\rightarrow} 0 \text{ as } n\rightarrow +\infty$ which implies 
$$\sqrt{n}(J(\widehat{\textbf{p}}_n,\textbf{q})-J(\textbf{p},\textbf{q}))\stackrel{\mathcal{D} }{ \rightsquigarrow} \mathcal{N}\left(0, V_{1,p} \right),\ \ \text{as}\ \ n\rightarrow+\infty.$$

\noindent This ends the proof of \eqref{thJ12n1}.

\bigskip \noindent The result \eqref{thJ12n2} is obtained by a symmetry argument by swapping the role of $ \textbf{p} $ and $\textbf{q}.$

\bigskip \noindent Now, it remains to prove Formula (\ref{thJ22n1}) of the theorem. Let us use bi-variate Taylor-Lagrange-Cauchy formula to get, 
\begin{eqnarray*}  
J(\widehat{\textbf{p}}_n,\widehat{\textbf{q}}_m)-J(\textbf{p},\textbf{q})
&=&\sum_{j\in D} \Delta_{p_n}^{c_j}\phi_{1}^{(1)}(p_j,q_j) + \sum_{j\in D}  \Delta_{q_m}^{c_j}\phi _{2}^{(1)}(p_j,q_j) \\
&&\frac{1}{2} \sum_{j\in D} \biggr((\Delta_{p_n}^{c_j})^{2}\phi^{(2)}_{1}+\Delta_{p_n}^{c_j}  \Delta_{q_m}^{c_j}\phi^{(2)}_{1,2}+ (\Delta_{q_m}^{c_j})^{2}\phi^{(2)}_{2}\biggr)\biggr(u_n^{c_j},v_m^{c_j} \biggr).   
\end{eqnarray*}

where $$
(u_n^{c_j},v_n^{c_j})=(p_j+\theta \Delta_{p_n}^{c_j}, \ q_j+\theta_j \Delta_{q_m}^{c_j}).
$$

\bigskip Thus we get  

\begin{eqnarray*}
J(\widehat{\textbf{p}}_n,\widehat{\textbf{q}}_m)-J(\textbf{p},\textbf{q}) &=&\frac{1}{\sqrt{n}}N_{n}(\textbf{p})+\frac{1}{\sqrt{m}}N_{m}(\textbf{q})+R_{n,m},
\end{eqnarray*}

\noindent where 
$$
N_{n}(\textbf{p})=\sum_{j\in D}\sqrt{p_j}\delta_n(p_j) \phi_{1}^{(1)}(p_{j},q_{j})\stackrel{\mathcal{D} }{\rightsquigarrow} \sum_{j\in D} \phi_{1}^{(1)}(p_{j},q_{j})Z_{p_{j}},\ \ \text{as}\ \ n\rightarrow+\infty
$$

$$
N_{m}(\textbf{q})=\sum_{j\in D}\sqrt{q_j}\delta_m(q_j)\phi _{2}^{(1)}(p_j,q_j) \stackrel{\mathcal{D} }{\rightsquigarrow} \sum_{j\in D} \phi_{2}^{(1)}(p_{j},q_{j})Z_{q_{j}},\ \ \text{as}\ \ m\rightarrow+\infty
$$

\noindent and $R_{n,m}$ is given by
$$
\frac{1}{2} \sum_{j\in D} \biggr((\Delta_{p_n}^{c_j})^{2}\phi^{(2)}_{1}+\Delta_{p_n}^{c_j}  \Delta_{q_m}^{c_j}\phi^{(2)}_{1,2}+ (\Delta_{q_m}^{c_j})^{2}\phi^{(2)}_{2}\biggr)\biggr(u_n^{c_j},v_m^{c_j} \biggr).
$$

\noindent First, we have that $N_{n}(\textbf{p})$ and $N_{m}(\textbf{p})$ are independents and hence

$$ 
N_{n}(\textbf{p})\stackrel{\mathcal{L} }{\sim} \mathcal{N}\left(0,V_{1,p} \right)\ \ \text{and} \ \ N_{m}(\textbf{q})\stackrel{\mathcal{L} }{\sim} \mathcal{N}\left(0,V_{2,q} \right).
$$

\noindent Therefore

\begin{eqnarray*}
\frac{1}{\sqrt{n}} \sum_{j\in D}\sqrt{p_j} \delta_n(p_j) \phi_{1}^{(1)}(p_{j},q_{j})+\frac{1}{\sqrt{m}} \sum_{j\in D}\sqrt{q_j}\delta_m(q_j)\phi_{2}^{(1)}(p_{j},q_{j})&=& \mathcal{N}\left(0,\frac{V_{1,p}}{n}+\frac{V_{2,q}}{m}\right)\\
	 &&\ \ \ \ +\ \  o_{\mathbb{P}}\left(\frac{1}{\sqrt{n}}\right)+o_{\mathbb{P}}\left(\frac{1}{\sqrt{m}}\right).
\end{eqnarray*}

\bigskip \noindent Thus

\begin{eqnarray*}
J(\widehat{\textbf{p}}_n,\widehat{\textbf{q}}_m)-J(\textbf{p},\textbf{q})&=&\mathcal{N}\left(0,\frac{V_{1,p}}{n}+\frac{V_{2,q}}{m}\right) +\ \  o_{\mathbb{P}}\left(\frac{1}{\sqrt{n}}\right)+o_{\mathbb{P}}\left(\frac{1}{\sqrt{m}}\right)+R_{n,m}.
\end{eqnarray*}

\noindent Next, we have 
\begin{eqnarray*}
\frac{1}{\sqrt{\frac{V_{1,p}}{n}+\frac{V_{2,q}}{m}}} \left( J(\widehat{\textbf{p}}_n,\widehat{\textbf{q}}_m)-J(\textbf{p},\textbf{q})\right)
=\mathcal{N}\left(0,1\right) &+& o_{\mathbb{P}}\left(\frac{1}{\sqrt{n}} \frac{1}{\sqrt{\frac{V_{1,p}}{n}+\frac{V_{2,q}}{m}}}\right)\\
	 &+& o_{\mathbb{P}}\left(\frac{1}{\sqrt{m}} \frac{1}{\sqrt{\frac{V_{1,p}}{n}+\frac{V_{2,q}}{m}}}\right)\\
&+& \ \ 
\frac{1}{\sqrt{\frac{V_{1,p}}{n}+\frac{V_{2,q}}{m}}}R_{n,m}.
\end{eqnarray*}

\noindent  That leads to  
\begin{eqnarray*}
\sqrt{\frac{nm}{mV_{1,p}+nV_{2,q}}}\left( J(\widehat{\textbf{p}}_n,\widehat{\textbf{q}}_m)-J(\textbf{p},\textbf{q})\right)&=&\mathcal{N}\left(0,1\right) +o_{\mathbb{P}}(1)+\sqrt{\frac{nm}{mV_{1,p}+nV_{2,q}}}R_{n,m},
\end{eqnarray*}

\noindent since $m/(m V_{1,p}+n V_{2,q})$ and $m/(n V_{1,p}+n V_{2,q})$ are bounded, and then
\begin{eqnarray*}
o_{\mathbb{P}}\left(\frac{1}{\sqrt{n}} \frac{1}{\sqrt{ \frac{V_{1,p}}{n}+\frac{V_{2,q}}{m}}}\right)&=&o_{\mathbb{P}}\left( \sqrt{\frac{m}{
 mV_{1,p}+nV_{2,q}}}\right)=o_{\mathbb{P}}(1)\\
&& and \\
o_{\mathbb{P}}\left(\frac{1}{\sqrt{m}} \frac{1}{\sqrt{ \frac{V_{1,p}}{n}+\frac{V_{2,q}}{m}}}\right)&=&o_{\mathbb{P}}\left( \sqrt{\frac{n}{
 mV_{1,p}+nV_{2,q}}}\right)=o_{\mathbb{P}}(1).
\end{eqnarray*}

\bigskip \noindent It remains to  prove  that $\left\vert 
\sqrt{\frac{nm}{mV_{1,p}+nV_{2,q}} }R_{n,m}\right\vert =o_{\mathbb{P}}(1).$ But we have by the continuity assumptions on $\phi$ and on its partial derivatives and by the uniform converges of $\Delta_{p_n}^{c_j}$ and $ \Delta_{q_m}^{c_j}$ to zero, that

\begin{eqnarray*}
&&\left\vert \sqrt{\frac{nm}{mV_{1,p}+nV_{2,q}} }R_{n,m}\right\vert \leq\\
&&\frac{1}{2} \left(\sqrt{n}a_{n}^{2} (\sum_{j\in D} \phi^{(2)}_{1}(p_j,q_j) +o(1))  \right) \left( \sqrt{\frac{m}{ mV_{1,p}+nV_{2,q}}}\right)\\
&+&\frac{1}{2} \left(\sqrt{m}b_{m}^{2} (\sum_{j\in D} \phi^{(2)}_{2}(p_j,q_j) +o(1))  \right) \left( \sqrt{\frac{n}{ mV_{1,p}+nV_{2,q}}}\right)\\
&+&\frac{1}{2} \left(\sqrt{n}a_{m}b_{m} (\sum_{j\in D} \phi^{(2)}_{2}(p_j,q_j) +o(1))  \right) \left( \sqrt{\frac{n}{ mV_{1,p}+nV_{2,q}}}\right)
\end{eqnarray*}

\bigskip \noindent As previously, we have $\sqrt{n}a_{n}^{2}=o_{\mathbb{P}}(1)$, $\sqrt{m}b_{m}^{2}=o_{\mathbb{P}}(1)$ and $\sqrt{n}a_{m}b_{m}=o_{\mathbb{P}}(1)$.

\bigskip \noindent From there, the conclusion is immediate.$\blacksquare$\\

\bigskip \noindent Proofs of Theorem \ref{thJs22} and Theorem \ref{thJs22d} are obtained by writing 
\begin{eqnarray*}
J^{(s)}(\widehat{\textbf{p}}_n,\textbf{q})-J^{(s)}(\textbf{p},\textbf{q})
&=&\frac{1}{2}\left(J(\widehat{\textbf{p}}_n,\textbf{q})-J(\textbf{p},\textbf{q})\right)+\frac{1}{2}\left(J(\textbf{q},\widehat{\textbf{p}}_n)-J(\textbf{q},\textbf{p})\right)
\\
 J^{(s)}(\textbf{p},\widehat{\textbf{q}}_m)-J^{(s)}(\textbf{p},\textbf{q}))&=&\frac{1}{2}\left(J(\textbf{p},\widehat{\textbf{q}}_m)-J(\textbf{p},\textbf{q})\right)+\frac{1}{2}\left(J(\widehat{ \textbf{q}}_m,\textbf{p})-J(\textbf{q},\textbf{p})\right)\\
  J^{(s)}(\widehat{\textbf{p}}_n,\widehat{\textbf{q}}_m)-J^{(s)}(\textbf{p},\textbf{q}))&=&\frac{1}{2}\left(J(\widehat{ \textbf{p}}_n,\widehat{\textbf{q}}_m)-J(\textbf{p},\textbf{q})\right)+\frac{1}{2}\left(J(\widehat{ \textbf{q}}_m,\widehat{\textbf{p}}_n)-J(\textbf{q},\textbf{p})\right) 
\end{eqnarray*}and using Theorems \ref{thJ12} and \ref{thJ22}.

\section{Simulations}
\label{simulation}

To assess the performance of ours estimators, we present a simulation study on a finite sample. In our
simulation study, we consider tree outcomes for the randoms variables $X$ and $Y$, $c_1,c_2,c_3$ with respectives \textit{p.m.f.}probabilities $p_1,
\ p_2,\, p_3$ and $q_1,\ q_2, \ q_3$. \\

\noindent Our aim is to compare the performance of the divergences measures  estimators as well as their symetrized forms with one or two samples 
when sample sizes increase. \\

\noindent Suppose that  $p_1=0.4,
\ p_2=0.25,\, p_3=0.35 $ and $q_1=0.27,\ q_2=0.32, \ q_3=0.41$.\\ 

\noindent With these above values of 
$p_j$ and $q_j$, $j=1,2,3$, the true and the symetrized form of our interest divergence measures become
\begin{eqnarray*}
&& \mathcal{D}_{T,0.99}(\textbf{p},\textbf{q})\approx 0.03969,\ \ \mathcal{D}_{T,0.99}^{(s)}(\textbf{p},\textbf{q}) \approx  0.03854\\
&&\mathcal{D}_{R,0.99}(\textbf{p},\textbf{q})\approx 0.03970,\ \ \mathcal{D}_{R,0.99}^{(s)}(\textbf{p},\textbf{q})\approx 0.03854\\
&& \mathcal{D}_{KL}(\textbf{p},\textbf{q})\approx 0.04012,\ \ \text{and}\ \ \mathcal{D}_{KL}^{(s)}(\textbf{p},\textbf{q}) \approx 0.03893.
\end{eqnarray*}

\bigskip \noindent 
\noindent In each figure, left panels represent the plots of divergence measure estimator, built from sample sizes of $n=100,200,\cdots,30000$, and the true divergence measure (represented by horizontal black line). The
middle panels show the histograms of the  data and where the red line represents the plots of the theoretical normal distribution calculated
from the same mean and the same standard deviation of the data. The right panels concern the Q-Q plot of the data which display the observed values against normally 
distributed data (represented by the red line). We see that the  underlying 
distribution of the data is normal since the points fall along a straight line.\\ 
 
\begin{figure}[H]
\centering
\includegraphics[scale=0.3]{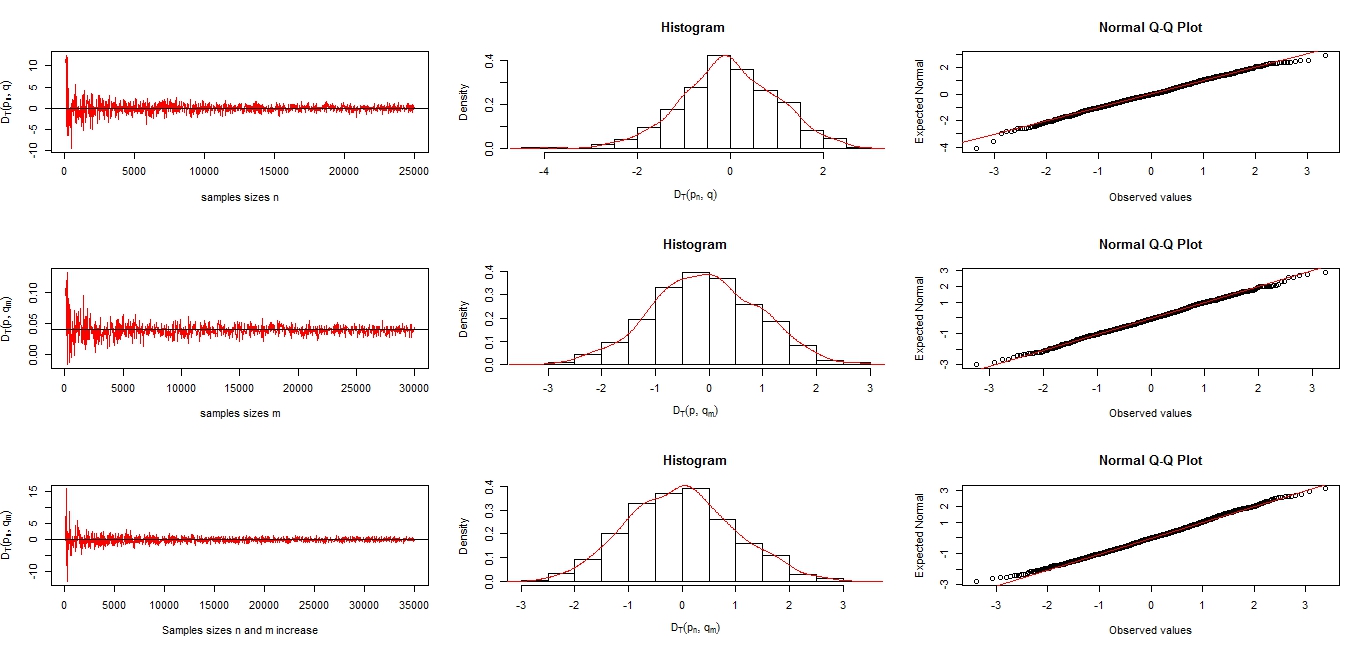} 
\caption{Plots when samples sizes increase, histograms and normal Q-Q plots of $\mathcal{D}_{T,\alpha}(p_n,q)$,   $\mathcal{D}_{T,\alpha}(p,q_m)$, and $\mathcal{D}_{T,\alpha}(p_n,q_m)$ ($\alpha=0.99$)
 versus $\mathcal{N}(0,1)$.
}\label{tsalpq}

\end{figure}

\begin{figure}[H]
\centering
\includegraphics[scale=0.3]{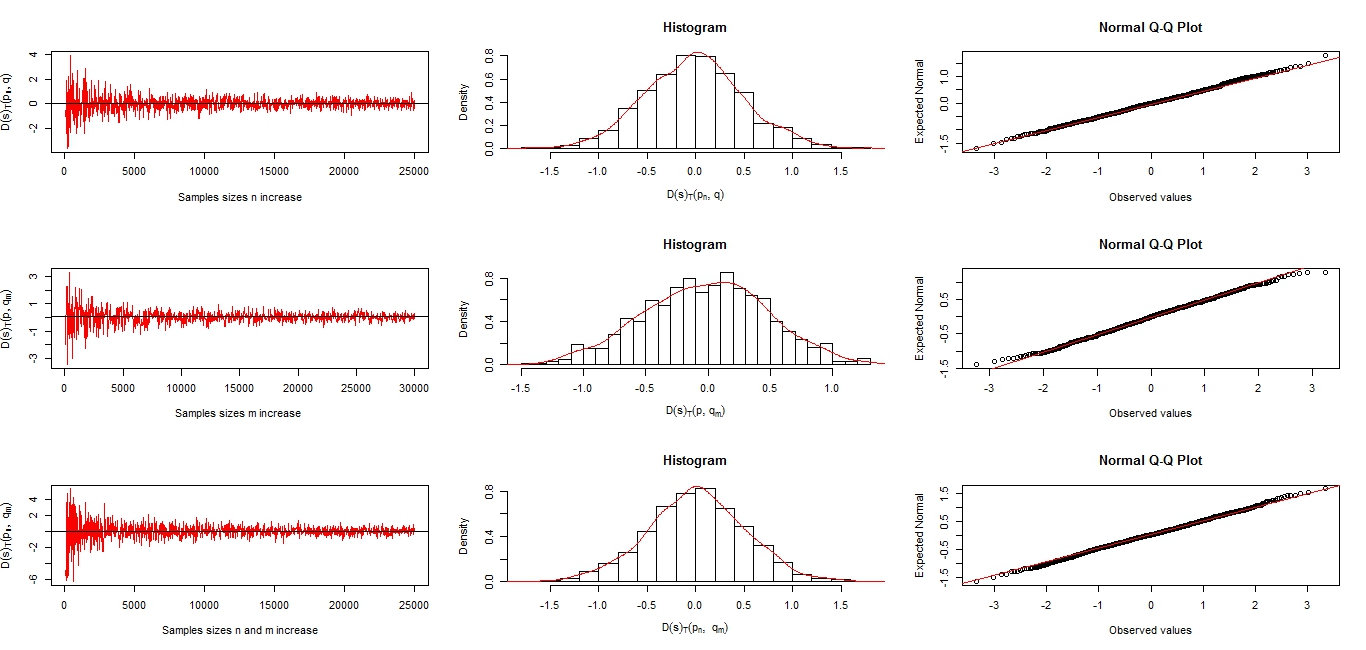} 
\caption{Plots when samples sizes increase, histograms and normal Q-Q plots of  $\mathcal{D}_{T,\alpha}^{(s)}(p_n,q)$,   $\mathcal{D}_{T,\alpha}^{(s)}(p,q_m)$, and $\mathcal{D}_{T,\alpha}^{(s)}(p_n,q_m)$ ($\alpha=0.99$) versus $\mathcal{N}(0,1)$.
}\label{tssympalpq}
\end{figure}

\begin{figure}[H]
\centering
\includegraphics[scale=0.3]{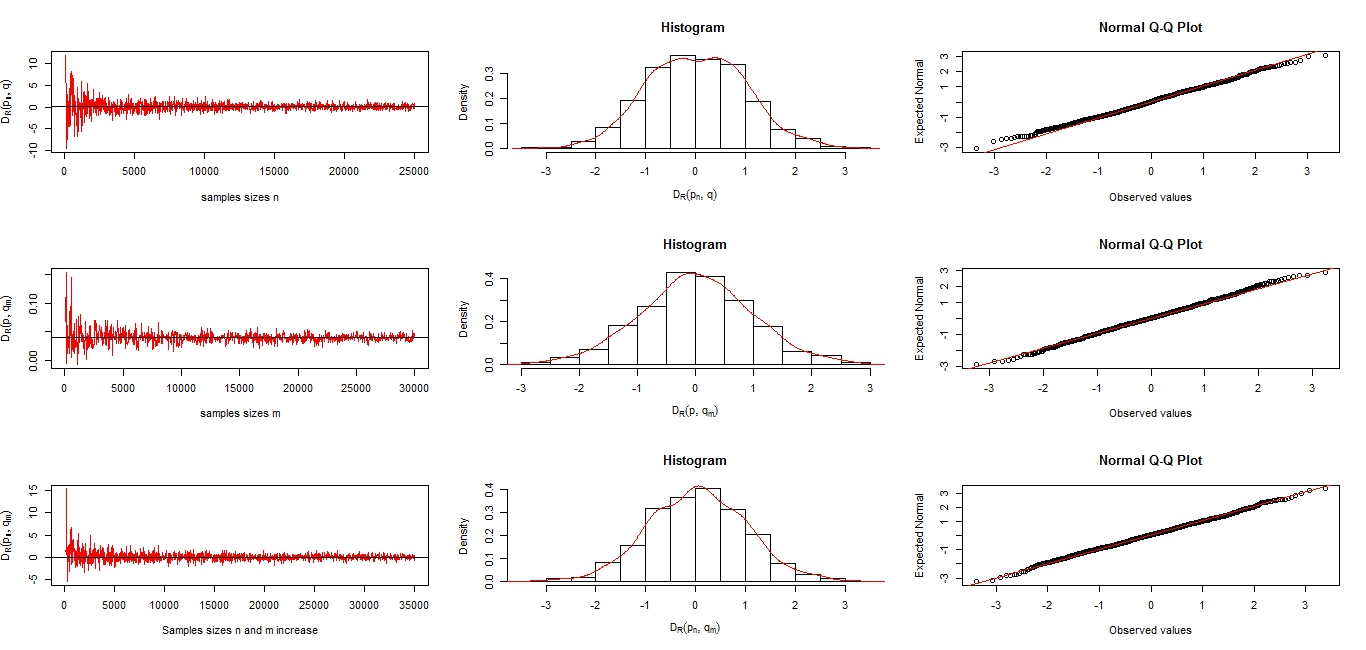} 
\caption{Plots when samples sizes increase, histograms and normal Q-Q plots of $\mathcal{D}_{R,\alpha}(p_n,q)$,   $\mathcal{D}_{R,\alpha}(p,q_m)$, and $\mathcal{D}_{R,\alpha}(p_n,q_m)$ ($\alpha=0.99$) versus $\mathcal{N}(0,1)$.
}\label{ralpq}
\end{figure}

\begin{figure}[H]
\centering
\includegraphics[scale=0.3]{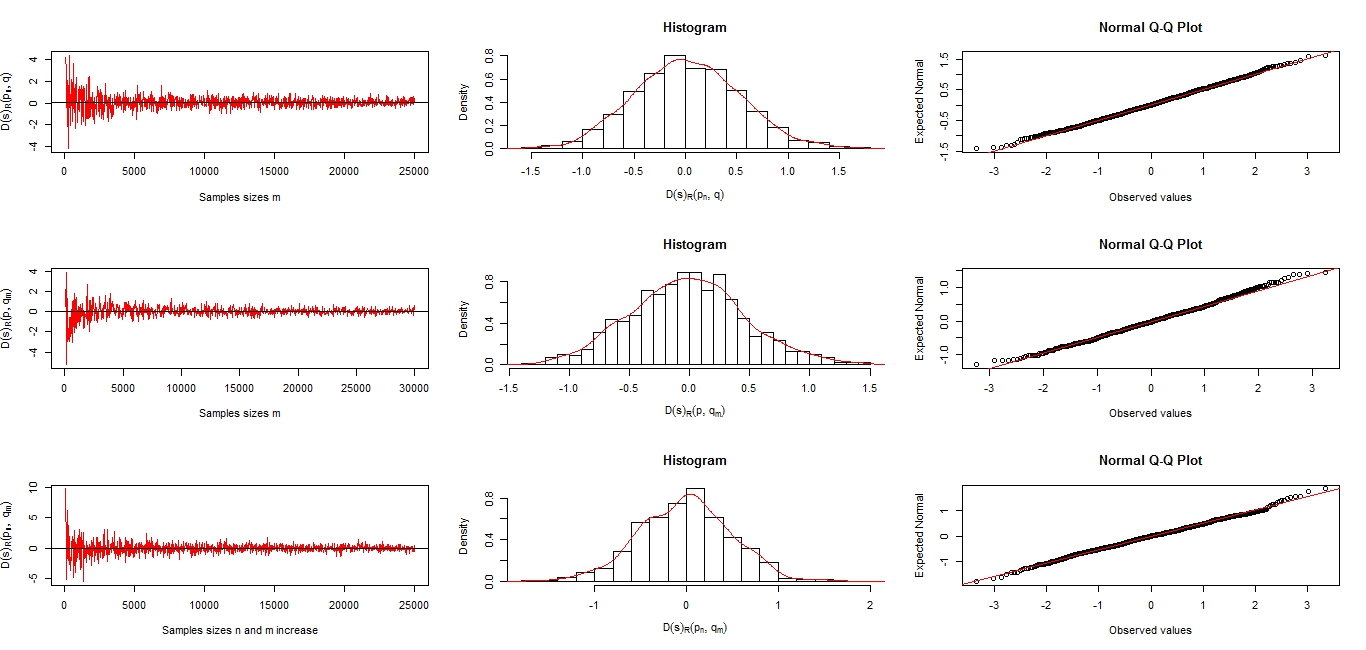} 
\caption{Plots  when samples sizes increase, histograms and normal Q-Q plots of  $\mathcal{D}_{R,\alpha}^{(s)}(p_n,q)$,   $\mathcal{D}_{R,\alpha}^{(s)}(p,q_m)$, and $\mathcal{D}_{R,\alpha}^{(s)}(p_n,q_m)$ ($\alpha=0.99$) versus $\mathcal{N}(0,1)$
}\label{rsympalpq}
\end{figure}

\begin{figure}[H]
\centering
\includegraphics[scale=0.3]{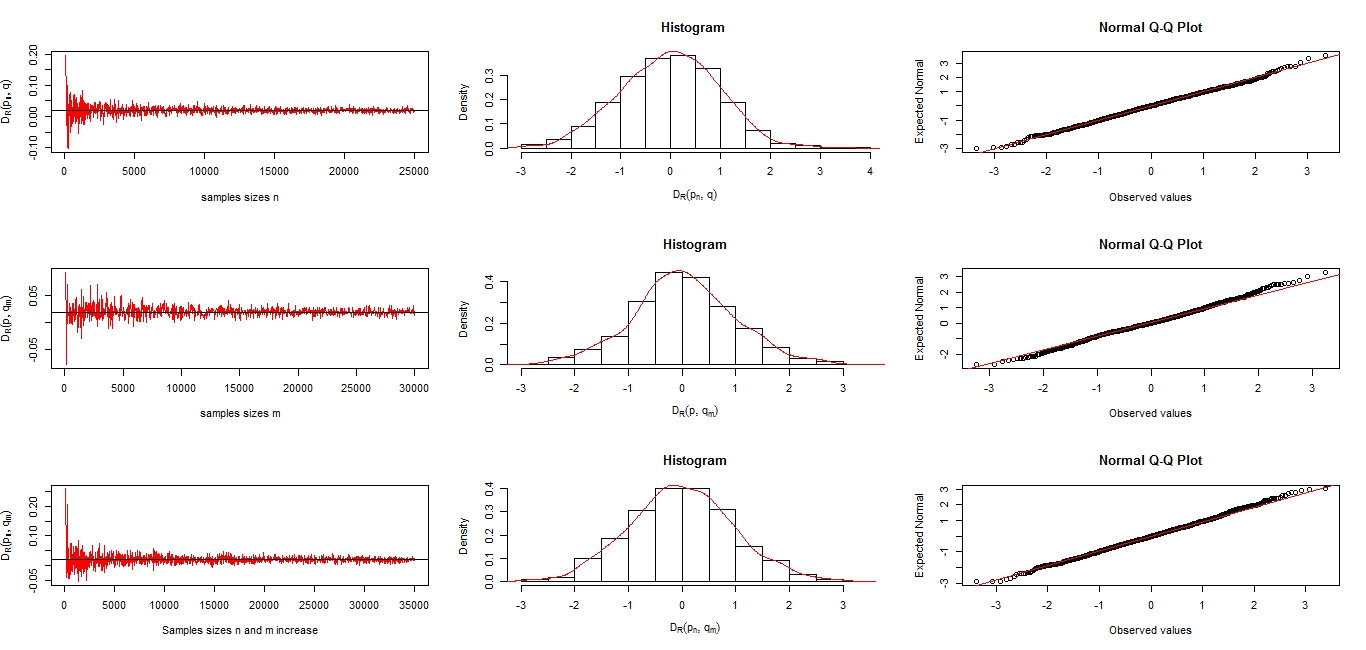} \caption{Plots when samples sizes increase, histogram and normal Q-Q plots of $\mathcal{D}_{R,\alpha}(p_n,q)$,   $\mathcal{D}_{R,\alpha}(p,q_m)$, and $\mathcal{D}_{R,\alpha}(p_n,q_m)$ ($\alpha=0.5$) versus $\mathcal{N}(0,1)$.
}\label{rsymalpq05}
\end{figure}

\begin{figure}[H]
\centering
\includegraphics[scale=0.3]{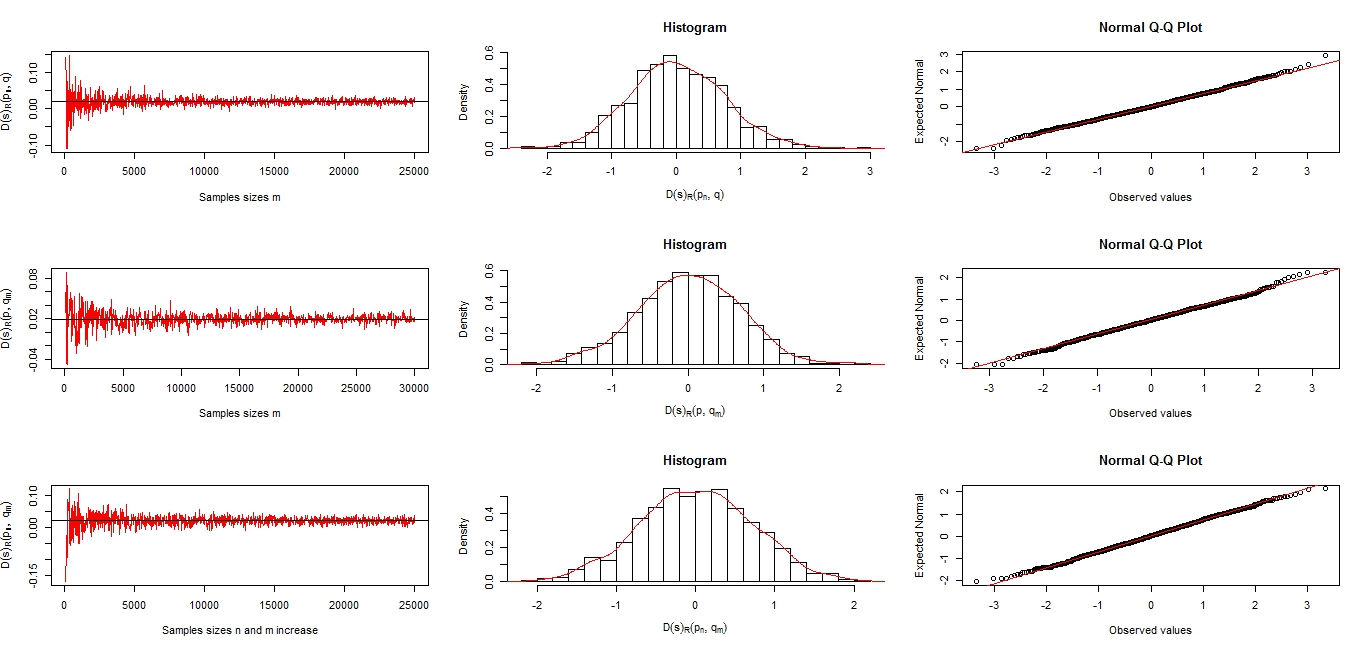} 
\caption{Plots when samples sizes increase, histograms and normal Q-Q plots of $\mathcal{D}_{R,\alpha}^{(s)}(p_n,q)$,   $\mathcal{D}_{R,\alpha}^{(s)}(p,q_m)$, and $\mathcal{D}_{R,\alpha}^{(s)}(p_n,q_m)$ ($\alpha=0.5$) versus $\mathcal{N}(0,1)$.
}\label{ralpq05}
\end{figure}

\begin{figure}[H]
\centering
\includegraphics[scale=0.3]{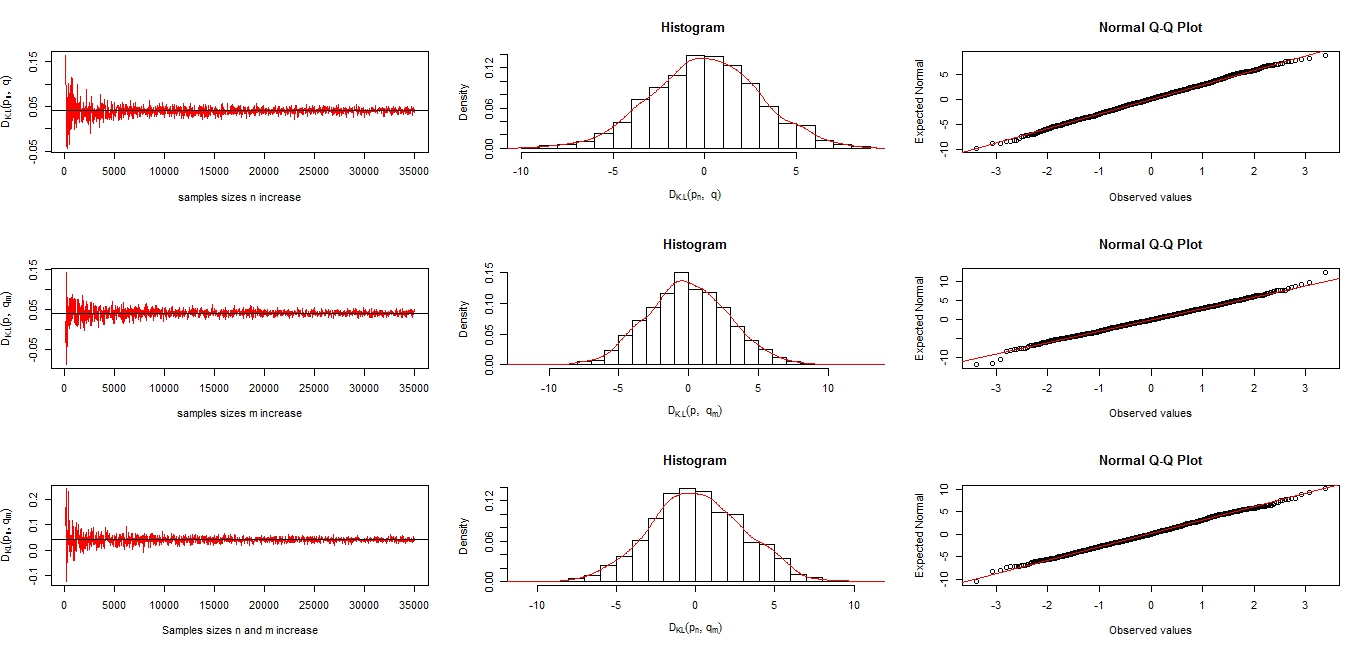}
\caption{Plots when samples sizes increase, histograms and normal Q-Q plots of $\mathcal{D}_{KL}(p_n,q)$,   $\mathcal{D}_{KL}(p,q_m)$, and $\mathcal{D}_{KL}(p_n,q_m)$ versus $\mathcal{N}(0,1)$.
}\label{pq}
\end{figure}

\begin{figure}[H]
\centering
\includegraphics[scale=0.3]{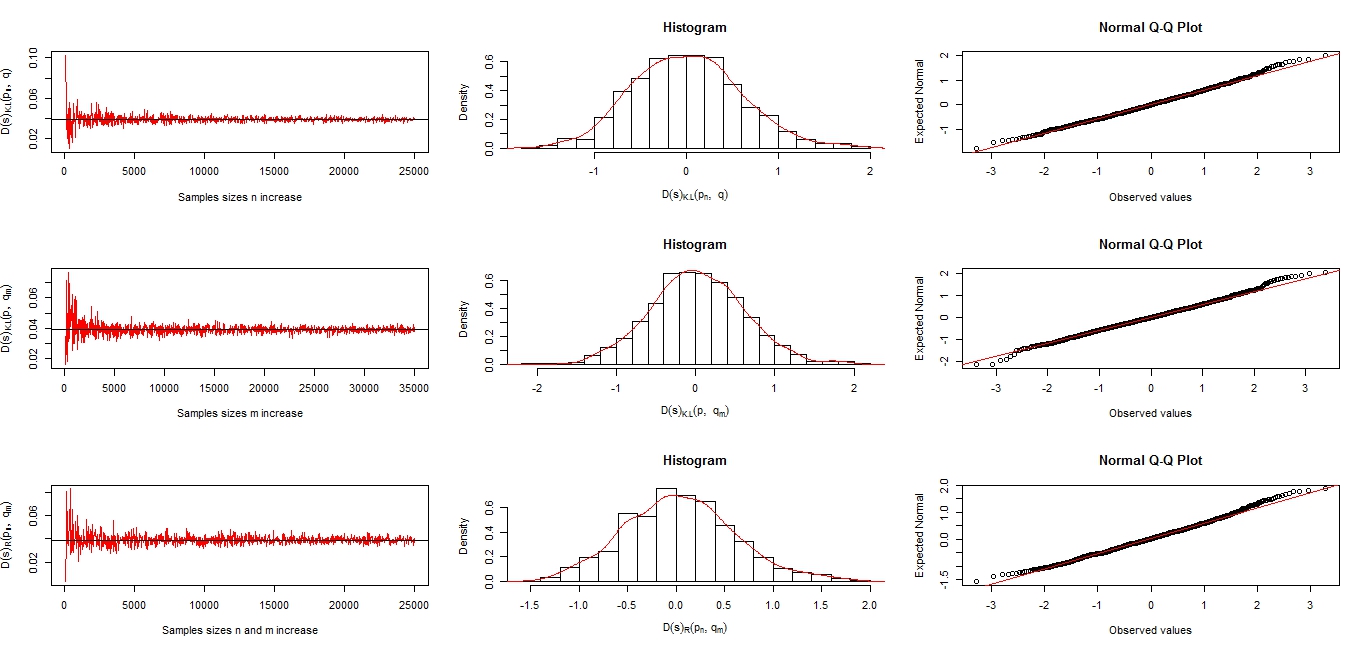} 
\caption{Plots when samples sizes increase, histograms and normal Q-Q plots of $\mathcal{D}_{KL}^{(s)}(p_n,q)$,   $\mathcal{D}_{KL}^{(s)}(p,q_m)$, and $\mathcal{D}_{KL}^{(s)}(p_n,q_m)$ versus $\mathcal{N}(0,1)$.
}\label{spq}

\end{figure}

  \noindent As seen in each plots in \textsc{figures} \ref{tsalpq},  \ref{tssympalpq}, \ref{ralpq},  \ref{rsympalpq}, \ref{rsymalpq05},  \ref{ralpq05}, \ref{pq}, and \ref{spq} our method performs well in showing the consistency of the divergence measures estimators and the asymptotic normality through their
 quantile-normal graphs.




\section{Conclusion}\label{conclus}
\noindent This paper joins a growing body of literature on estimating divergence measures in
the discrete case and on finite sets. We adopted the plug-in method
and we derived almost sure rates of convergence and asymptotic normality of the
most common divergence measures in one sample, two samples as well as symetrical
form of divergence measures, all this, by means of the functional $\phi-$divergence
measure.

\end{document}